\documentclass[fleqn,10pt]{SelfArx} 
\usepackage[english]{babel} 
\usepackage{lipsum} 
\usepackage{cite}
\usepackage{alltt}
\usepackage{moreverb}
\usepackage[]{algorithm,algcompatible,lipsum}
\usepackage[]{algpseudocode}
\usepackage{listings}
\usepackage{multirow}
\usepackage{flushend}

\newcommand{\norma}[1]{\left \| #1 \right \|}

\newcommand{\dpara}[2]{\frac{\partial #1}{\partial #2}}
\newcommand{\dparb}[3]{\frac{\partial #1}{\partial #2 
\partial #3}}

\newcommand{\dparc}[4]{\frac{\partial #1}{\partial #2 \partial #3
\partial #4}}

\newcommand{\dpard}[5]{\frac{\partial #1}{\partial #2 \partial #3
\partial #4 \partial #5}}

\newcommand{\gbf}[1] {\mbox{\boldmath${#1}$\unboldmath}}

\definecolor{color1}{RGB}{0,0,90} 
\definecolor{color2}{RGB}{0,20,20} 

\usepackage{xurl} 
\usepackage{hyperref} 

\hypersetup{
	hidelinks,
	colorlinks,
	breaklinks=true,
	urlcolor=color1,
	citecolor=color1,
	linkcolor=color1,
	bookmarksopen=false,
	pdftitle={Title},
	pdfauthor={Author},
}

\Archive{} 

\PaperTitle{A dual number formulation to efficiently compute higher order
directional derivatives} 

\Authors{R. Pe\'on-Escalante\textsuperscript{1,2}, 
K.~B. Cant\'un-Avila\textsuperscript{1},
O. Carvente\textsuperscript{1},
A. Espinosa-Romero\textsuperscript{2},
F. Pe\~nu\~nuri\textsuperscript{1}*
} 
\affiliation{\textsuperscript{1}\textit{Facultad de Ingenier\'ia, Universidad Aut\'onoma de Yucat\'an, A.P. 150, Cordemex, M\'erida, Yucat\'an, M\'exico.}} 

\affiliation{\textsuperscript{2}\textit{Facultad de Matem\'aticas, Universidad Aut\'onoma de 
Yucat\'an, Anillo Perif\'erico Norte, Tablaje Cat. $13615$, 
Colonia Chuburn\'a Hidalgo Inn, M\'erida Yucat\'an, M\'exico.}} 

\affiliation{*\textbf{Corresponding author}: fracisco.pa@correo.uady.mx} %

\Keywords{Directional derivatives, dual numbers, automatic 
differentiation, multilinear forms, kinematics quantities,
 finite differences.} 
\newcommand{\keywordname}{Keywords} 

\Abstract{This contribution proposes a new formulation to efficiently compute 
directional derivatives of order one to fourth. The formulation is based 
on automatic differentiation implemented with dual numbers. Directional 
derivatives are particular cases of symmetric multilinear forms; 
therefore, using their symmetric properties and their coordinate 
representation, we implement functions to calculate mixed partial 
derivatives. Moreover, with directional derivatives, we deduce concise 
formulas for the velocity, acceleration, jerk, and jounce/snap vectors.  
The utility of our formulation is proved with three examples. The first 
example presents a comparison against the forward mode of finite 
differences to compute the fourth-order directional derivative of a 
scalar function. To this end, we have coded the finite differences 
method to calculate partial derivatives until the fourth order, to any 
order of approximation. The second example presents efficient 
computations of the velocity, acceleration, jerk, and jounce/snap. 
Finally, the third example is related to the computation of some partial 
derivatives. The implemented code of the proposed formulation and the 
finite differences method is proportioned as additional material to this 
article.}

\begin{document}

\maketitle 

\tableofcontents 

\thispagestyle{empty} 

\section{Introduction}\label{sec1} 

The computation of derivatives is a subject of great relevance in many 
areas of  science and engineering. Derivatives are required by many 
numerical methods as well as to optimize functions, compute 
velocities and accelerations, compute sensitivities, among many other 
applications \cite{Margossian2019, Atilim2018,Penunuri2020,
Chandrasekhar2021,Dvurechensky2021}. The finite differences method (FD) 
\cite{boole2009}, has been the most used method to compute any sort of 
derivatives. However, this method is subject to truncation and 
cancellation errors which, in general, are extremely difficult to 
overcome. Although at first order, the complex step approximation method 
\cite{Martins2003} solves the cancellation problem, this persists for 
higher-order derivatives. Despite these drawbacks, no one can deny that
FD is of paramount importance in both theoretical and practical studies. 
One of its main  advantages is that FD can be directly applied to 
differentiate numerical data \cite{KHAN2003,SOD19781,Martins2013,
BALAKRISHNA2021}. Automatic differentiation (AD) is another method to 
compute derivatives numerically which, unlike FD, is not subject to  
truncation or cancellation errors \cite{Griewank1989, Bischof1992,
CHINCHALKAR1994197, GAWA2003, Neidinger2010, Guo2021, 
OBERBICHLER2021113817}. The accuracy of the resulting derivatives is 
only related to the computer precision determined by its architecture. It could 
be said that AD cannot deal with the differentiation of numerical data 
because the function to generate them is unknown \cite{KHAN2003}; 
however, AD can be used if the unknown function is constructed, for 
instance, using spline interpolation of the data as was done in 
\cite{CANTUNAVILA2021104086}.

The first studies on AD can be traced back around the decade of sixties 
\cite{Juedes1991,Griewank2008}, although just in recent years, there 
have been a large number of studies using it in many areas of the 
sciences and engineering, including but not limited to robotics, Monte 
Carlo data analysis, mechanisms, machine learning, heat transfer, 
differential-algebraic equations, solid mechanics, quantum mechanics, 
and chemical physics \cite{Giftthaler2017,RAMOS201919,Peon2020,
Bolte2020,NIU2021,ESTEVEZSCHWARZ2021, Vigliotti2021,Fakher2022,
Kasim2022}. There are essentially two ways of implementing AD: the 
forward and the backward mode, some bibliography on the implementation 
of AD can be seen in \cite{10.1007/3-540-28438-9_28}. The forward mode 
can be implemented with dual numbers  \cite{Cheng1994,Jeffrey2012,
Penunuri2013,Wenbin2013,Mendoza2015,Mex2015}. Most studies implementing 
AD with dual numbers deal with first or second-order derivatives. 
Moreover, the second-order case generally addresses single variable 
functions, with some exceptions in which mixed partial derivatives are 
considered  \cite{TANAKA2015,Penunuri2020}. One of the contributions of 
our study is to present and code formulas to compute mixed partial 
derivatives up to the fourth order, enhancing even more the potential 
applications of AD. The presented formulation could be extended to 
arbitrary order; nevertheless, in physics and engineering, higher orders 
than four are rarely necessary, if any.

The dual number implementation of AD is exceptionally efficient when 
computing directional derivatives in general rather than the partial 
derivatives by themselves. The first and second-order directional 
derivative of a function $f:\mathbb{R}^m \to \mathbb{R}$ along the 
vector $\mathbf{v}$, evaluated at the point $\mathbf{q}$  are the most 
common scalar directional derivatives, which involves products of the 
gradient and Hessian operators with a vector;  they are given by:
$\mathbf{D}_{\mathbf{v}}f(\mathbf{q})$ $=$ $\nabla f(\mathbf{q})\,
\mathbf{v}$ and $\mathbf{D}^2_{\mathbf{v}}f(\mathbf{q})$ $=$ 
$\mathbf{v}^\text{T} \mathbf{H}f(\mathbf{q}) \mathbf{v},$ respectively. 
It is worth noting that these quantities are frequently needed instead 
of merely gradients or Hessians \cite{Moller1993,Pearlmutter1994,
Hicken2014,Doug2018,Song2022}. So, instead of calculating such operators 
and their products with vector $\mathbf{v}$ the dual number 
implementation directly computes these directional derivatives with a 
single evaluation of the function  $f$  in a dual variable 
\cite{Penunuri2020}.

Directional derivatives can also be vector quantities. The product of 
the Jacobian matrix with a vector is the first-order vector directional 
derivative. Higher order vector directional derivatives arise (as we 
will show), for instance, when computing the jerk \cite{Schot1978} and 
the jounce/snap \cite{Eager2016}--the third and fourth order time 
derivative of the position vector of a point particle--, quantities 
(especially the jerk) of relevance in engineering \cite{Eager2016,
Figliolini2019,Fang2020}.
The main contribution of this study is to construct a formulation based 
on dual numbers to calculate scalar and vector directional derivatives 
up to the fourth order along several vectors, for instance, 
$\mathbf{v}^\text{T} \mathbf{H} f(\mathbf{q}) \mathbf{w}$  instead of 
$\mathbf{v}^\text{T} \mathbf{H} f(\mathbf{q}) \mathbf{v}$. The case of 
directional derivatives along a single vector can be directly
addressed by taking the dual parts of a dual function. The case of
directional derivatives along several vectors is more complicated, and 
we are not aware of a methodology to address this case. We have found 
that the problem can be solved using the properties of multilinear forms
\cite{Hassani1999,Carroll2004}, from which directional derivatives are a 
particular case. In order to implement the proposed formulation to 
compute higher-order directional derivatives, we develop some 
general-purpose Fortran modules and Matlab classes to compute these 
kinds of symmetric multilinear forms. Our implementation can cope with
complex and real cases. This is worth mentioning since few works support
AD for complex numbers  \cite{Guo2021}. The Fortran and Matlab 
implementations are independent; the former allows us to conduct heavy 
computations if necessary. The efficiency of the proposed formulation to
compute directional derivatives is contrasted with a traditional 
approach using finite differences. For this reason, we implement the 
forward mode of the finite differences (FMFD) method in Fortran to 
compute partial derivatives until the fourth order to any order of 
approximation. This is the third contribution of our study since the 
computational implementation to calculate higher order partial 
derivatives is difficult to find,  especially to any order of 
approximation, and with the possibility of using single, double, and 
quadruple precision variables. With the proposed formulation, all the 
directional derivatives, from order one to fourth, are computed with a 
single evaluation of the given function, overloaded to accept dual 
variables.

As an additional contribution to this study, there are included some 
examples. The first shows the efficiency of the proposed formulation 
compared to FD or a symbolic approach to computing a fourth-order 
directional derivative. In the second example, formulas for all the 
kinematics quantities up to the fourth order, namely, velocity, 
acceleration, jerk, and jounce/snap, in terms of directional derivatives
are deduced. Also, we present an alternative to the directional 
derivative approach to compute such kinematics quantities by using a 
Taylor expansion of the vector function of generalized coordinates. The 
last example presents the computation of some partial derivatives using 
dual numbers and the FMFD method. The code of our implementation and the 
presented examples in  this study are available at \cite{r_peon2022M}.

The paper is organized as follows. Section 
\ref{sec2} introduces the main definitions relating to multilinear 
forms and how higher-order directional derivatives are 
particular cases of them. Section \ref{sec3} summarizes the forward mode 
of the finite differences method, there we present how this method can 
be coded to any order of approximation. Section \ref{sec4} describes how 
dual numbers can be used to efficiently compute higher order directional 
derivatives. This section also explains how to compute partial 
derivatives until the fourth order. Some examples are presented in 
Section \ref{Ejemplos} and finally the conclusions are given in Section 
\ref{conclusions}.

\section{Multilinear forms and higher order directional derivatives}
\label{sec2}
This section presents some concepts related to multilinear forms which 
are essentials in the deduction of formulas to compute directional
derivatives using dual numbers.
\subsection{Multilinear forms}
Let $r$ be a natural number and $X_1$, $X_2$, $\dots$, $X_r$, and $Y$ 
vector spaces. A multilinear form, multilinear map or multilinear 
transformation, is a map 
\begin{equation}
\mathbf{G}: X_1\times X_2 \dots \times X_r \to Y,
\label{eq:1}
\end{equation}
with the following properties.
\begin{align}
\mathbf{G}&(\mathbf{x}_1, \dots,\mathbf{x}_i + \mathbf{x'}_i,
\dots,\mathbf{x}_r) = \mathbf{G}(\mathbf{x}_1, \dots,
\mathbf{x}_i,\dots,\mathbf{x}_r) + \nonumber\\
&+\mathbf{G}(\mathbf{x}_1, \dots,
\mathbf{x'}_i,\dots,\mathbf{x}_r)
\end{align}
\begin{align}
\mathbf{G}(\mathbf{x}_1, \dots,\lambda \, \mathbf{x}_i,\dots,
\mathbf{x}_r) &= \lambda \, \mathbf{G}(\mathbf{x}_1, \dots,
\mathbf{x}_i,\dots,\mathbf{x}_r),
\end{align}
for all $\mathbf{x}_i\in X_i$, $\mathbf{x'}_i \in X_i$ for 
$i \in \{1,\dots,r\}$, and $\lambda \in \mathbb{R}$.

In this study we are interested in the particular case where 
the vector spaces $X_i$ are $\mathbb{R}^m$ and $Y$ is $\mathbb{R}^n$. In
the case of $Y=\mathbb{R}$ the multilinear form receive the name of 
tensor \cite{Hassani1999,Carroll2004}. We first present the discussion 
for this particular case as can be directly generalized to the 
case $Y=\mathbb{R}^n$. Indeed, choosing the standard basis 
$\mathbf{e}_\alpha$ of 
$\mathbb{R}^n$
we have
%
\begin{align}
\mathbf{G}(\mathbf{x}_1, \dots,\mathbf{x}_i + \mathbf{x'}_i,
\dots,\mathbf{x}_r) &= \sum_{\alpha=1}^n g_\alpha(\mathbf{x}_1, \dots,
\mathbf{x}_i + \mathbf{x'}_i,\nonumber\\
&\dots,\mathbf{x}_r)\,\mathbf{e}_\alpha,
\end{align}
%
where $g_\alpha$ is a multilinear map $g:\mathcal{R}^{rm} \to 
\mathbb{R}$ (a tensor), with
\begin{align}
\mathcal{R}^{rm} := \underbrace{\mathbb{R}^m \times \dots 
\times\mathbb{R}^m}_{\mbox{$r$ times}}.
\end{align}
Therefore, the multilinear form we are interested in is:
\begin{align}
&g:\mathcal{R}^{rm} \to \mathbb{R},~ \text{with the properties} 
\nonumber \\ 
&g(\mathbf{x}_1, \dots,\mathbf{x}_i + \mathbf{x'}_i,
\dots,\mathbf{x}_r) = g(\mathbf{x}_1, \dots,\mathbf{x}_i,\dots,
\mathbf{x}_r) + \nonumber\\
&~~+g(\mathbf{x}_1, \dots,\mathbf{x'}_i,\dots,\mathbf{x}_r)
\nonumber \\
&g(\mathbf{x}_1, \dots,\lambda\,\mathbf{x}_i,
\dots,\mathbf{x}_r) = \lambda\,g(\mathbf{x}_1, \dots,\mathbf{x}_i,\dots,
\mathbf{x}_r).
\end{align}

The number $r$ defines the rank of a multilinear form, in this study
we take $r\in\{1,2,3,4\}$. This is because we are interested in 
computing partial derivatives until the fourth order, hence, all the
functions are assumed to be at least of class C$^4$ (four times 
differentiable in all its variables). 

A multilinear form is symmetric if 
\begin{align}
g(\mathbf{x}_1, \dots, \mathbf{x}_r) = g(\pi(\mathbf{x}_1, \dots,
\mathbf{x}_r)),
\end{align}
where $\pi(\mathbf{x}_1, \dots,\mathbf{x}_r)$ is any
permutation of $\{\mathbf{x}_1, \dots,\mathbf{x}_r\}$. 
Choosing a particular basis for $\mathbb{R}^m$ we can construct the 
coordinate representation of a multilinear transformation.
Thus, the coordinate representation for $r=4$ is
\begin{align}\label{matrep}
A_{ijkl}=g(\mathbf{e}_{i},\mathbf{e}_{j},\mathbf{e}_{k},\mathbf{e}_{l}),
\end{align}
where we have used the standard basis of $\mathbb{R}^m$.
For the cases $r=1,2,3$ we have
\begin{align}\label{matrep1}
A_{i}&=g(\mathbf{e}_{i})\\
\label{matrep2}
A_{ij}&=g(\mathbf{e}_{i},\mathbf{e}_{j})\\
\label{matrep3}
A_{ijk}&=g(\mathbf{e}_{i},\mathbf{e}_{j},\mathbf{e}_{k}).
\end{align}

With Eq. (\ref{matrep}) the multilinear map of rank 4 is rewritten as:
\begin{align}
g(\mathbf{x}, \mathbf{y}, \mathbf{z}, \mathbf{w}) = 
\sum_{i,j,k,l=1}^m{A_{ijkl}x_i y_j z_k w_l},
\end{align}
which will be written as
\begin{align} \label{R4MF}
g(\mathbf{x}, \mathbf{y}, \mathbf{z}, \mathbf{w}) = 
A_{ijkl}x_i y_j z_k w_l,
\end{align}
where $v_i$ is the $i$-$th$ component of vector $\mathbf{v}$ and we have
used the Einstein summation convention, which will be used in this study 
with the exception of \ref{A1}. The cases $r=1,2,3$ are written below.

\begin{align}\label{R1MF}
g(\mathbf{x}) &= A_{i}x_i \\
\label{R2MF}
g(\mathbf{x}, \mathbf{y}) &= A_{ij}x_i y_j \\
\label{R3MF}
g(\mathbf{x}, \mathbf{y}, \mathbf{z}) &= A_{ijk}x_i y_j z_k 
\end{align}

\subsection{Higher order directional derivatives}
The first order directional derivative of a function 
$f:\mathbb{R}^m \to \mathbb{R}$ along the vector 
$\mathbf{x} \in \mathbb{R} ^m$ (in this study $\mathbf{x}$ is not 
necessary a unit vector) evaluated 
at the point $\mathbf{q} \in \mathbb{R} ^m$ is 
defined as
\begin{equation}
D^1_{\mathbf{x}}f(\mathbf{q}) = D_\mathbf{x} f(\mathbf{q}) = 
\lim_{h\to 0} \frac{f(\mathbf{q} + h\,\mathbf{x}) - f(\mathbf{q})}{h}
\label{DirDer1Lim}
\end{equation}
Since we are assuming all the functions of class--at least--C$^4$, 
expanding the numerator of this limit in a Taylor series, we can prove 
that
\begin{equation}
 D_\mathbf{x} f(\mathbf{q}) = \nabla f(\mathbf{q}) \,
\mathbf{x} = \dpara{f(\mathbf{q})}{q_i}\,x_i.
\label{DirDer1Grad}
\end{equation}
Higher-order directional derivatives can be recursively defined as 
\cite{Marban1969,Hiriart1984,Seeger1988}
\begin{equation}
D^r_{\mathbf{x}}f(\mathbf{q}) = 
D_{\mathbf{x}}[D^{r-1}_{\mathbf{x}}f(\mathbf{q})]=
\nabla[D^{r-1}_{\mathbf{x}}f(\mathbf{q})] \,\mathbf{x}.
\end{equation}
From this, after applying Eq. (\ref{DirDer1Lim}) or 
Eq. (\ref{DirDer1Grad}), the second, third and fourth order directional 
derivatives are given by
\begin{align}
\label{dx2DD}
D^2_{\mathbf{x}}f(\mathbf{q}) &=
\dparb{^2 f(\mathbf{q})}{q_i}{q_j} x_i x_j \\ 
\label{dx3DD}
D^3_{\mathbf{x}}f(\mathbf{q}) &=
\dparc{^3 f(\mathbf{q})}{q_i}{q_j}{q_k} x_i x_j x_k \\
\label{dx4DD}
D^4_{\mathbf{x}}f(\mathbf{q}) &=
\dpard{^4 f(\mathbf{q})}{q_i}{q_j}{q_k}{q_l} x_i x_j x_k x_l.
\end{align}
The above higher-order directional derivatives along the vector 
$\mathbf{x}$ can be generalized to directional derivatives along 
several vectors. For instance, the 
second order directional derivative along two different vectors  
is constructed as follows. 
Applying Eq. (\ref{DirDer1Lim}) to $d(\mathbf{q}) = 
D_\mathbf{y} f(\mathbf{q}) = \nabla f(\mathbf{q}) \,
\mathbf{y}$, we have
\begin{align}\label{gd2A}
D_{\mathbf{x}}[D_{\mathbf{y}}f(\mathbf{q})] = 
D^2_{\mathbf{x,y}}f(\mathbf{q})=
\lim_{h\to 0} 
\frac{d(\mathbf{q} + h\,\mathbf{x}) - d(\mathbf{q})}{h} \\
\label{gd2B}
D^2_{\mathbf{x,y}}f(\mathbf{q}) = \lim_{h\to 0} 
\frac{[\nabla f(\mathbf{q} + h\,\mathbf{x})] \mathbf{y} - 
\nabla f(\mathbf{q})\,\mathbf{y}}{h}.
\end{align}
Now, expanding the numerator of Eq. (\ref{gd2B}) in a Taylor series it 
follows that
\begin{equation}
D^2_{\mathbf{x,y}}f(\mathbf{q})= 
\dparb{^2 f(\mathbf{q})}{q_i}{q_j} x_i y_j.
\label{gd2Def2}
\end{equation}
Comparing with Eq. (\ref{R2MF}) we see that Eq. (\ref{gd2Def2}) is a
symmetric multilinear form of rank 2. 
For convenience, we change the notation accordingly to a multilinear 
form,
therefore, we define
\begin{align}
\label{Dtodfq1}
D_{\mathbf{x}}f(\mathbf{q})&= d_{1f\mathbf{q}}(\mathbf{x}) \\
\label{Dtodfq2}
D^2_{\mathbf{x,y}}f(\mathbf{q})&= d_{2f\mathbf{q}}(\mathbf{x,y})\\
\label{Dtodfq3}
D^3_{\mathbf{x,y,z}}f(\mathbf{q})&= d_{3f\mathbf{q}}(\mathbf{x,y,z})\\
\label{Dtodfq4}
D^4_{\mathbf{x,y,z,w}}f(\mathbf{q})&=d_{4f\mathbf{q}}(\mathbf{x,y,z,w}).
\end{align}
Formulas for the third and fourth order directional derivatives along 
\mbox{several} vectors can be obtained in similar way to the second 
order case. These directional derivatives (including the first and 
second order) are given by:
\begin{align}\label{dmf1}
d_{1f\mathbf{q}}(\mathbf{x}) &=  \dpara{f\,(\mathbf{q})}{q_i}\,x_i \\
\label{dmf2}
d_{2f\mathbf{q}}(\mathbf{x,y}) &=  \dparb{^2 f (\mathbf{q})}{q_i}{q_j}\,
x_i\,y_j \\
\label{dmf3}
d_{3f\mathbf{q}}(\mathbf{x,y,z}) &=  \dparc{^3 f(\mathbf{q})}{q_i}{q_j}{q_k}
\,x_i\,y_j\,z_k \\
\label{dmf4}
d_{4f\mathbf{q}}(\mathbf{x,y,z,w}) &=  \dpard{^4 f(\mathbf{q})}{q_i}{q_j}{q_k}{q_l}
\,x_i\,y_j\,z_k\,w_l.
\end{align}
Since the functions are at least of class C$^4$, all these directional 
derivatives are symmetric multilinear forms involving partial 
derivatives of ranks 1, 2, 3 and 4. In physics and engineering ranks 
higher than 4 are rarely needed. This is the main reason we consider 
$r\in\{1,2,3,4\}$. 
From now on we will use the phrase: higher order directional derivative, 
for both cases; when such a directional derivative is computed along a 
single vector or when it is computed along several vectors. From the 
arguments of the involved multilinear form it will be clear what case 
are we talking about. Also, to lighten the notation, the dependence on 
$f$ and $\mathbf{q}$ will be obviated in the left hand side of 
Eqs. (\ref{dmf1}--\ref{dmf4}). For instance, the 
rank 4 multilinear form involving partial derivatives (fourth order 
directional derivative) will be written as
\begin{align}
d_4(\mathbf{x,y,z,w}) &=  \dpard{^4 f (\mathbf{q})}{q_i}{q_j}{q_k}{q_l}\,
x_i\,y_j\,z_k\,w_l.
\label{d4mlffor}
\end{align}

A traditional approach, using finite differences (FD), to compute Eqs. 
(\ref{dmf1}--\ref{dmf4}) is in general 
inefficient. Moreover the inherent subtraction and cancellation error of 
FD, could cause serious precision lost on the computed values. An 
alternative is to use dual numbers to implement automatic 
\mbox{differentiation}. The advantage of implementing AD with dual 
numbers is that the aforementioned equations are computed efficiently, 
without truncation or cancellation errors. 

\section{Forward mode of finite differences to compute derivatives}
\label{sec3}
Before we develop the dual number approach to compute higher order 
directional derivatives, we briefly present the finite differences 
method. Finite differences is a powerful tool to compute the derivatives 
of a function. Unfortunately, to obtain an accurate value of a 
derivative, the computational program may requires to define variables 
of high precision. This produces inefficient code to compute derivatives. 
Moreover, the computation of the directional derivatives need to 
calculate all the involved derivatives in their sums. Nevertheless,
the utility of FD is undeniable. This section summarizes the forward 
mode of the finite differences method (FMFD) to compute partial 
derivatives. Most of the equations here presented, although scattered, 
can be found elsewhere in the literature, nevertheless, we recommend 
\cite{boole2009}. On the other hand the computational implementation to 
compute higher order partial derivatives to any order of approximation
is difficult to find. For the single variable case, however, we 
recommend \cite{Hassan2012}.

\subsection{The forward delta operator}
The forward delta operator $\gbf{\Delta}_h$ is defined as
\begin{align}\label{Dhdef}
\gbf{\Delta}_h = \mathbf{T}_h -\mathbf{I},
\end{align}
with
\begin{align}
\mathbf{T}_h\,f(q) &= f(q+h)\\
\mathbf{I}f(q) &= f(q),
\end{align}
being $f$ a function $f:\mathbb{R}\to \mathbb{R}$.

From the Taylor series of $f(q)$ and denoting the derivative operator as
$D$ it can be seen that
\begin{align}
f(q+h)-f(q) &= h f'(q) + \frac{h^2}{2!}f''(q) +\dots \\
\gbf{\Delta}_h f(q) &= \left[ h D + 
\frac{h^2}{2!}D^2 + \dots \right] f(q),
\end{align}
from these equations we see that
\begin{align}
\gbf{\Delta}_h = e^{hD}-\mathbf{I},
\end{align}
and thus the derivative operator can be written as
\begin{align}\label{Doplog}
D = \frac{1}{h} \log(\gbf{\Delta}_h + \mathbf{I})
\end{align}
or
\begin{align}\label{hderap}
h\,D =\gbf{\Delta}_h - \frac{1}{2}\gbf{\Delta}_h^2+
\frac{1}{3}\gbf{\Delta}_h^3-\frac{1}{4}\gbf{\Delta}_h^4 + \dots
\end{align}

Using Eq. (\ref{Dhdef}) it is not difficult to prove that
\begin{align}\label{Dop}
\gbf{\Delta}_h ^p f(q) = \sum_{k=0}^p (-1)^{p-k} \,
\frac{p!}{(p-k)!\,k!}\, f(q + hk).
\end{align}

\subsection{Derivatives of function of one variable}
From Eq. (\ref{hderap}) we can approximate the first order 
derivative of $f$ to any order of approximation $n$.
Defining the operator
\begin{align}\label{hfppdef}
\gbf{\delta}_h^n &= \left[\gbf{\Delta}_h - \frac{1}{2}
\gbf{\Delta}_h^2+\frac{1}{3}\gbf{\Delta}_h^3  + \dots + 
\frac{(-1)^{n+1}}{n}\gbf{\Delta}_h^n\right]
\end{align}
we have
\begin{align}\label{hfpp}
\frac{1}{h}\gbf{\delta}_h^n f(q) &= f'(q) + O(h^{n})~\text{or}\\
 f'(q) &=\frac{1}{h}\gbf{\delta}_h^n f(q) + O(h^{n})
\end{align}
and then
\begin{align}\label{hfppm1}
\gbf{\delta}_h^n f(q) = h f'(q) + O(h^{n+1})
\end{align}

Higher order derivatives can be computed with successive applications of
the operator $\gbf{\delta}_h^n$, for instance an approximation to order 
8 of $f''(q)$ could be obtained as
\begin{align}
&\gbf{\delta}_h^{4} \cdot \gbf{\delta}_h^{5} f(q) =  h^2 f''(q) + 
O\left(h^{10}\right)\\
&\frac{1}{h^2}\gbf{\delta}_h^{4} \cdot \gbf{\delta}_h^{5} f(q) = f''(q)+ 
O\left(h^{8}\right)
\end{align}
In general, we have
\begin{align}\label{hpderpn1p}
&\gbf{\delta}_h^{n_1} \cdot \gbf{\delta}_h^{n_2} \dots \cdot 
\gbf{\delta}_h^{n_r} f(q) = h^r\,f^{(r)}(q) \;+ \nonumber\\
&~~+\;O\left(h^{n_1+n_2+\dots + n_r+1}\right)\\
\label{derpn1p}
&\frac{1}{h^r}\gbf{\delta}_h^{n_1} \cdot \gbf{\delta}_h^{n_2} \dots \cdot 
\gbf{\delta}_h^{n_r} f(q) = f^{(r)}(q) + \nonumber\\
&~~+\;O\left(h^{n_1+n_2+\dots + n_r+1-r}\right).
\end{align} 

\subsection{Partial derivatives of function of several variables}
Eqs. (\ref{hpderpn1p}, \ref{derpn1p}) can be directly generalized to 
compute partial derivatives of a function of several variables 
$f:\mathbb{R}^m \to \mathbb{R}$. This generalization is
\begin{align}\label{dpars}
&\frac{\gbf{\delta}_{h_{ir}}^{n_1} \dots \cdot \gbf{\delta}_{h_{i1}}^
{n_r} f(\mathbf{q})}{h_{i1} \dots h_{ir}}  = \frac{\partial ^r}
{\partial q_{i1}\dots\partial q_{ir}}f(\mathbf{q}) \;+ \nonumber\\
&~~+\;O\left(\norma{\mathbf{h}}^{n_1+\dots + n_r+1-r}\right)
\end{align}
with $\mathbf{h}=[h_{i1},\dots,h_{ir}]$ and 
\begin{align}
\gbf{\Delta}_{h_j} = f(\mathbf{q} + h_j \mathbf{e}_j)-f(\mathbf{q}),~~
j \in \{i1,i2,\dots, ir\}.
\end{align} 

Eq. (\ref{dpars}) can be implemented nesting the required functions as 
follow. After coding Eq. (\ref{Dop}),  Eqs. (\ref{hfppdef}, \ref{hfpp}) 
follow directly. Let  {\tt deriv} the function which implements the
action of the $\gbf{\delta}_{h}^{n}$ operator on a function, explicitly,
\begin{align}
{\tt deriv}(f, \mathbf{q}_0, h, n) =\frac{1}{h} \gbf{\delta}_{h}^{n}
f(\mathbf{q}_0),
\end{align}
then, the partial derivatives $(\partial /\partial q_i) f(\mathbf{q}_0)$
can be calculated with the function of Alg. \ref{fdd1}. 
%

\begin{algorithm}[htb]
\caption{Partial derivatives, $\partial f/\partial q_i$ evaluated at 
$\mathbf{q}_0$ with an approximation order $n$.}
\label{fdd1}
\begin{algorithmic}[1]
\Function{fr=df1}{$f,\mathbf{q}_0,i, h, n$}
 \State fr $\gets$ \verb+deriv+(faux,$\mathbf{q}_0(i)$,$h$,$n$)
\Function{$f_0$=faux}{$q$}
\State \hspace{0.0cm}$\mathbf{v}$ $\gets$ $\mathbf{q}_0$
\State \hspace{0.0cm}$\mathbf{v}(i)$ $\gets$ $q$
\State \hspace{0.0cm}$f_0$ $\gets$ f($\mathbf{v}$)
\EndFunction
\EndFunction
\end{algorithmic}
\end{algorithm}
%

This procedure can be applied to compute higher order derivatives. For 
instance the second and third order partial derivatives can be computed 
using Algorithms \ref{fdd2} and \ref{fdd3} respectively. In these 
algorithms we have used $h_{i1}=h_{i2}=h_{i3}=h$.

\begin{algorithm}[htb]
\caption{Second order partial derivatives, $\partial^2 f/\partial q_j
\partial q_i$ evaluated at $\mathbf{q}_0$.}
\label{fdd2}
\begin{algorithmic}[1]
\Function{fr=df2}{$f,\mathbf{q}_0,i,j,h,n_1,n_2$}
\State fr $\gets$ \verb+deriv+(faux,$\mathbf{q}_0(j)$,$h$,$n_1$)
\Function{$f_1$=faux}{$q$}
\State \hspace{0.0cm}$\mathbf{v}$ $\gets$ $\mathbf{q}_0$
\State \hspace{0.0cm}$\mathbf{v}(j)$ $\gets$ $q$
\State \hspace{0.0cm}$f_1$ $\gets$ df1($f,\mathbf{v},i,h,n_2$)
\EndFunction
\EndFunction
\end{algorithmic}
\end{algorithm}

\begin{algorithm}[htb]
\caption{Third order partial derivatives, $\partial^3 f/\partial q_k 
\partial q_j \partial q_i$ evaluated at $\mathbf{q}_0$. }
\label{fdd3}
\begin{algorithmic}[1]
\Function{fr=df3}{$f,\mathbf{q}_0,i,j,k,h,n_1,n_2,n_3$}
\State fr $\gets$ \verb+deriv+(faux,$\mathbf{q}_0(k)$,$h$,$n_1$)
\Function{$f_2$=faux}{$q$}
\State \hspace{0.0cm}$\mathbf{v}$ $\gets$ $\mathbf{q}_0$
\State \hspace{0.0cm}$\mathbf{v}(k)$ $\gets$ $q$
\State \hspace{0.0cm}$f_2$ $\gets$ df2($f,\mathbf{v},i,j,h,n_2,n_3$)
\EndFunction
\EndFunction
\end{algorithmic}
\end{algorithm}

Observe, from Eq. (\ref{dpars}), that for a given $r$ and a desired 
order of approximation $o_a$, the equation
\begin{equation}
n_1+n_2+\dots + n_r = o_a +r -1
\label{eq:nsum}
\end{equation}
must be satisfy. Therefore the set of numbers $\{n_1,+n_2,\dots,  n_r\}$
is not unique. For instance for $r=4$ and $o_a=8$ two possible sets of 
values are $\{8,1,1,1\}$ and $\{2,3,3,3\}$. In principle, for 
computations with infinite precision, both set of values would give the 
same result. Nevertheless, in our numerical experiments, the second set 
of values (when all the $n_i$ are the closest to each other and sorted 
in ascending order) gives a more accurate result. This set of values can 
be obtained using the function of algorithm \ref{setvn}. Our Fortran 
implemented code to compute partial derivatives, until fourth order (to 
any order of approximation) using FMFD can be found in \cite{r_peon2022M}. 
It is worthwhile to mention that, from a practical point of view, 
computing the derivatives to a very high order of approximation will 
cause errors of precision. The computational implementation allows 
variables of single precision (real32), double precision (real64) and 
quadruple precision (real128).

\begin{algorithm}
\caption{Function to obtain a possible set of values to satisfy a given 
approximation order $o_a$.}
\label{setvn}
\begin{algorithmic}[1]
\Function{$\mathbf{fr}$=setvn}{$r,o_a$}
\State res $\gets$  module($o_a + r - 1, r$)
\State nmin $\gets$  $(o_a + r - 1 -\text{res})/r$
\For {$k= 1$, $r-\text{res}$}
	\State $\mathbf{fr}$($k$)  $\gets$ nmin
\EndFor
\For {$k= r-\text{res}+1$, $r$}
	\State $\mathbf{fr}$($k$)  $\gets$ $\text{nmin}+1$
\EndFor
\EndFunction
\end{algorithmic}
\end{algorithm}
\section{Dual number to compute derivatives} \label{sec4}
\subsection{Dual numbers and derivatives of scalar functions of one 
variable}

In the late XIX century, Clifford \cite{Clif1873} introduces  a number 
of the form
\begin{align}\label{dual1N}
\hat{q} = a_0 + a_1\epsilon_1
\end{align}
with $a_0$ and $a_1$ real numbers and $\epsilon_1$ the dual entity with
the property \mbox{$\epsilon_1\cdot \epsilon_1=0$}. Dual numbers are not 
a field since a pure dual number does not have multiplicative inverse. 
Nevertheless, dual numbers form a commutative ring which can be used to 
compute derivatives. This can be seen expanding the Taylor series of an 
analytic function $f:\mathbb{R}\to \mathbb{R}$ (or $f:\mathbb{C}\to 
\mathbb{C}$ for the complex case), 
\begin{align}
f(q+h) &= f(q) + h f'(q) + \frac{h^2}{2} f''(q) + O(h^3)\\
\label{EQdual1}
f(q+\epsilon_1) &= f(q) + f'(q)\,\epsilon_1, \text{ since }\epsilon_1^2 
= 0.
\end{align}
Observe that the dual part of Eq. (\ref{EQdual1}) is the derivative of 
$f$.
The computation of higher order derivatives can be conducted defining
\begin{align} \label{dualn}
\hat{z} = \sum_{k=0}^n a_k \epsilon_k, 
\end{align}
$\text{with }\epsilon_0=1$ and $a_k$ real (or complex) numbers. Now, 
evaluating the Taylor expansion of $f(x)$ to order $n$ in the dual 
number $x+\epsilon_1$ we can see that, in order to obtain the 
derivatives of $f$, we must chose
\begin{equation}\label{eq:eik}
\epsilon_1^k  =
\begin{cases}
k!\,\epsilon_k & \mbox{if } k \leqslant n,
\vspace*{0.3cm}\\
0 & \mbox{if } k>n.
\end{cases}
\end{equation}
Then multiplying $\epsilon_i$ by $\epsilon_j$ 
and assuming $i+j\leqslant n$ (if $i+j > n$ the product is 0) we have
\begin{equation}
\epsilon_i \cdot \epsilon_j = \frac{\epsilon_1 ^i}{i!}\,
\frac{\epsilon_1 ^j}{j!} = \frac{\epsilon_1 ^{i+j}}{i!\,j!} = 
\frac{(i+j)!}{i!\,j!} \epsilon_{i+j}
\label{eq:eik2}
\end{equation}
thus
\begin{equation}\label{gentabmult}
\epsilon_i \cdot \epsilon_j =
\begin{cases}
0 & \mbox{if } i+j > n,
\vspace*{0.3cm}\\
\frac{(i+j)!}{i!\,j!} \epsilon_{i+j} & \mbox{otherwise. }
\end{cases}
\end{equation}
A different deduction of this equation can be seen in \cite{Kalos2021}.

Although from a theoretical point of view the generalization to compute 
derivatives of any order is straightforward, its computational 
implementation requires considerably effort. Most of the implementation 
of dual numbers in the literature are for orders 1 or 2. Since in this 
study we are interesting in computing generalized directional 
derivatives until order four, we have implemented  Eq. (\ref{dualn}) 
up to this order. 

\subsection{Dual numbers to compute multilinear forms involving partial 
derivatives} \label{sec:DNDD}
\subsubsection{Higher order directional derivatives: case 
$f:\mathbb{R}^m \to \mathbb{R}$. }
Since we have implemented Eq. (\ref{dualn}) for the complex case,
the below discussion also applies when the involved vector space is 
$\mathbb{C}^n$. However, for the sake of concreteness, the discussion 
is done talking about $\mathbb{R}$ instead of $\mathbb{C}$.

From the Taylor series expansion until fourth order of a function 
$f:\mathbb{R}^m\to \mathbb{R}$ we have
\begin{align}\label{TaylorS4f}
&f(\mathbf{q}+ \mathbf{h})= f(\mathbf{q}) + h_i \dpara{f(\mathbf{q})}
{q_i} + \frac{h_i h_j}{2!} \dparb{^2 f(\mathbf{q})}{q_i}{q_j} \;+ \nonumber\\
&~~+\;\frac{h_i h_j h_k}{3!} \dparc{^3 f(\mathbf{q})}{q_i}{q_j}{q_k} +\frac{h_i h_j h_k h_l}{4!} \dpard{^4 f(\mathbf{q})}
{q_i}{q_j}{q_k}{q_l} + \nonumber\\
&~~+\;R(\mathbf{q},\mathbf{h}), \end{align}
with $R(\mathbf{q},\mathbf{h})$
of the order $\norma{\mathbf{h}}^5$. Overloading the function to 
accept dual arguments, evaluating in the point 
$\mathbf{q}+\epsilon_1\mathbf{v}$ and from Eq. (\ref{gentabmult}) 
with $n=4$ we have
\begin{align}\label{FvecDual4} 
&f(\mathbf{q}+\epsilon_1\mathbf{v})=\epsilon_0 f(\mathbf{q}) + 
\epsilon_1 \dpara{f(\mathbf{q})}{q_i}v_i  +\epsilon_2\, 
\dparb{^2 f(\mathbf{q})}{q_i}{q_j} v_i v_j + \nonumber\\
&~~+\;\epsilon_3 
\dparc{^3 f(\mathbf{q})}{q_i}{q_j}{q_k} v_i v_j v_k + \epsilon_4\, \dpard{^4 f(\mathbf{q})}
{q_i}{q_j}{q_k}{q_l} v_i v_j v_k v_l,
\end{align}
which form Eqs. (\ref{dmf1}--\ref{dmf4}) can be 
rewritten as
\begin{align}\label{fqmasev}
&f(\mathbf{q}+\epsilon_1\mathbf{v})=\epsilon_0 f(\mathbf{q}) + 
\epsilon_1 d_1(\mathbf{v})  + \epsilon_2 d_2(\mathbf{v},\mathbf{v}) +
\nonumber \\
&~~+ \epsilon_3 d_{3}(\mathbf{v},\mathbf{v},\mathbf{v}) + 
\epsilon_4 d_{4}(\mathbf{v},\mathbf{v},\mathbf{v},\mathbf{v}).
\end{align}
Now, defining $f. \epsilon _k$ the $\epsilon_k$ part of $f$
we have
\begin{align}\label{d1v}
d_{1}(\mathbf{v})&=f(\mathbf{q}+\epsilon_1\mathbf{v}).\epsilon_1\\
\label{d2vv}
d_{2}(\mathbf{v},\mathbf{v})&=f(\mathbf{q}+\epsilon_1\mathbf{v}).
\epsilon_2\\
\label{d3vvv}
d_{3}(\mathbf{v},\mathbf{v},\mathbf{v})&=f(\mathbf{q}+
\epsilon_1\mathbf{v}).\epsilon_3\\
\label{d4vvvv}
d_{4}(\mathbf{v},\mathbf{v},\mathbf{v},\mathbf{v})&=f(\mathbf{q}+
\epsilon_1\mathbf{v}).\epsilon_4
\end{align}
All of these expressions are efficiently computed with only 
one evaluation of the function $f$. More importantly, there will not be 
truncation or cancellation errors when computing them. Observe that 
Eqs. (\ref{d2vv}--\ref{d4vvvv}) are particular cases of 
Eqs. (\ref{dmf2}--\ref{dmf4}). Nevertheless, using the 
properties of a symmetric multilinear form, these particular cases can 
be applied to compute the general case. This is achieved as follows.

Computing $d_2(\mathbf{x}+\mathbf{y},\mathbf{x}+\mathbf{y})$ we have
\begin{align}
&d_2(\mathbf{x}+\mathbf{y},\mathbf{x}+\mathbf{y}) = 
d_2(\mathbf{x},\mathbf{x}) + d_2(\mathbf{x},\mathbf{y}) + 
d_2(\mathbf{y},\mathbf{x}) + \nonumber \\
&~~+d_2(\mathbf{y},\mathbf{y})
\end{align}
since the involved multilinear forms are symmetric we have
\begin{align}
d_2(\mathbf{x},\mathbf{y}) = \frac{1}{2}\left[
d_2(\mathbf{x}+\mathbf{y},\mathbf{x}+\mathbf{y}) - 
d_2(\mathbf{x},\mathbf{x}) - d_2(\mathbf{y},\mathbf{y})\right]
\end{align}
The right hand side of this equation can be computed using 
Eq. (\ref{d2vv}). A similar procedure can be used to prove that 
\begin{align}\label{d3xyz}
&d_3(\mathbf{x},\mathbf{y},\mathbf{z}) = \frac{1}{2}\Big[
d_3(\mathbf{x}+\mathbf{y},\mathbf{x}+\mathbf{y},\mathbf{z}) -
d_3(\mathbf{x},\mathbf{x},\mathbf{z}) +\nonumber \\
&~~ -d_3(\mathbf{y},\mathbf{y},\mathbf{z})\Big],
\end{align}
with
\begin{align}\label{d3xxz}
&d_3(\mathbf{x},\mathbf{x},\mathbf{z})=\frac{1}{6}\Big[
d_3(\mathbf{z} + \mathbf{x},\mathbf{z} + \mathbf{x},\mathbf{z} +
\mathbf{x}) + \nonumber\\
&~~+d_3(\mathbf{z} - \mathbf{x},\mathbf{z} - 
\mathbf{x},\mathbf{z} - \mathbf{x}) - 2 d_3(\mathbf{z},\mathbf{z},
\mathbf{z})\Big],
\end{align}
which can be computed using Eq. (\ref{d3vvv}).

To compute $d_4(\mathbf{x},\mathbf{y},\mathbf{z},\mathbf{w})$, we 
first define 
\begin{align}
d_4xz^{+} = d_4(\mathbf{x}+\mathbf{z},\mathbf{x}+\mathbf{z},\mathbf{x}+
\mathbf{z},\mathbf{x}+\mathbf{z})\\
d_4xz^{-} = d_4(\mathbf{x}-\mathbf{z},\mathbf{x}-\mathbf{z},\mathbf{x}-
\mathbf{z},\mathbf{x}-\mathbf{z})
\end{align}
and note that
\begin{align}
&d_4xz^{+} + d_4xz^{-}= 12 d_4(\mathbf{x},\mathbf{x},\mathbf{z},
\mathbf{z})  + 2 d_4(\mathbf{x},\mathbf{x},\mathbf{x},\mathbf{x}) + \nonumber\\
&~~+2 d_4(\mathbf{z},\mathbf{z},\mathbf{z},\mathbf{z}),
\end{align}
thus
\begin{align}
&d_4(\mathbf{x},\mathbf{x},\mathbf{z},\mathbf{z}) = \frac{1}{12}\Big[
d_4xz^{+} + d_4xz^{-} - 2 d_4(\mathbf{x},\mathbf{x},\mathbf{x},
\mathbf{x}) +\nonumber\\
&~~ - 2 d_4(\mathbf{z},\mathbf{z},\mathbf{z},\mathbf{z})\Big],
\end{align}
Now computing $d_4(\mathbf{x},\mathbf{x},\mathbf{z}+\mathbf{w},
\mathbf{z}+\mathbf{w})$ it can be proved that
\begin{align}
&d_4(\mathbf{x},\mathbf{x},\mathbf{z},\mathbf{w})=\frac{1}{2}\Big[
d_4(\mathbf{x},\mathbf{x},\mathbf{z}+\mathbf{w},\mathbf{z}+\mathbf{w})+\nonumber\\
&~~-
d_4(\mathbf{x},\mathbf{x},\mathbf{z},\mathbf{z})-
d_4(\mathbf{x},\mathbf{x},\mathbf{w},\mathbf{w})\Big]
\end{align}
then computing $d_4(\mathbf{x}+\mathbf{y},\mathbf{x}+\mathbf{y},
\mathbf{z},\mathbf{w})$ we obtain
\begin{align}
&d_4(\mathbf{x},\mathbf{y},\mathbf{z},\mathbf{w}) = \frac{1}{2}\Big[
d_4(\mathbf{x}+\mathbf{y},\mathbf{x}+\mathbf{y},\mathbf{z},\mathbf{w})+\nonumber\\
&~~ - d_4(\mathbf{x},\mathbf{x},\mathbf{z},\mathbf{w})-
d_4(\mathbf{y},\mathbf{y},\mathbf{z},\mathbf{w})\Big],
\end{align}
which in the last instance, it can be computed using Eq. (\ref{d4vvvv}).

We have formulas to compute all the Eqs. (\ref{dmf1}--\ref{dmf4}). It is interesting to note that the product of 
the Hessian matrix with two vectors is given by
\begin{align}
\mathbf{v}^\text{T} \mathbf{H}f(\mathbf{q}) \mathbf{w} = d_2(\mathbf{v},
\mathbf{w}),
\end{align}
and then, the $k$-th component of the product of the Hessian matrix with 
a vector is given by
\begin{align}
\left[ \mathbf{H}f(\mathbf{q}) \mathbf{w}\right]_k = d_2(\mathbf{e}_k,
\mathbf{w}).
\end{align}
Moreover using the coordinate representation of a multilinear form we 
can compute all the mixed partial derivatives until fourth order. For 
instance,
\begin{align}\label{ejemplod3}
\dparc{^3 f}{q_i}{q_j}{q_k}=d_3(\mathbf{e}_i,\mathbf{e}_j,\mathbf{e}_k),
~\text{for}~i,j,k \in \{1,2 \dots,m\}.
\end{align}

Algorithm \ref{dirderr3} shows how to compute the directional derivative 
of order 3 along 3 different vectors--Eq. (\ref{d3xyz}).  Using this 
algorithm, Eq. (\ref{ejemplod3}) can be computed with algorithm 
\ref{parderd3}. There, the vector $\mathbf{indxv}$ stores the variable 
indexes $i,j,k$. The evaluating point is $\mathbf{q} \in \mathbb{R}^m$. 
The directional derivatives of orders 1, 2 and 4 can be coded in similar 
way.

\begin{algorithm}[htb]
\caption{Directional derivative of third order (multilinear form of 
rank 3).}
\label{dirderr3}
\begin{algorithmic}[1]
\Function{fr=d3mlf}{$f,\mathbf{x},\mathbf{y},\mathbf{z},\mathbf{q}$}
\State fr $\gets$ 0.5\Big[d3qxz($f,\mathbf{x} + \mathbf{y},\mathbf{z},
\mathbf{q}$) $-$ d3qxz($f,\mathbf{x},\mathbf{z},\mathbf{q}$) $-$ 
d3qxz($f,\mathbf{y},\mathbf{z},\mathbf{q}$)\Big]
\EndFunction

\Function{fr=d3qxz}{$f,\mathbf{x},\mathbf{z},\mathbf{q}$}
\State fr $\gets$ $\frac{\text{d3q}(f,\mathbf{z}+\mathbf{x},\mathbf{q})+
\text{d3q}(f,\mathbf{z}-\mathbf{x},\mathbf{q}) - 
2\,\text{d3q}(f,\mathbf{z},\mathbf{q})}{6}$
\EndFunction

\Function{fr=d3q}{$f,\mathbf{v},\mathbf{q}$}
\State fqd $\gets$ $f(\mathbf{q} + \epsilon_1\,\mathbf{v})$
\State fr $\gets$ fqd$.\epsilon_3$
\EndFunction
\end{algorithmic}
\end{algorithm}

\begin{algorithm}[htb]
\caption{Third order partial derivatives, $\partial^3 f/\partial q_k 
\partial q_j \partial q_i$. Coordinate representation of the
rank 3 multilinear form, associated to the third order directional 
derivative computed with dual numbers.}
\label{parderd3}
\begin{algorithmic}[1]
\Function{fr=df3}{$f,\mathbf{indxv},\mathbf{q}$}
\State $i$ $\gets$ $\mathbf{indxv}(1)$
\State $j$ $\gets$ $\mathbf{indxv}(2)$
\State $k$ $\gets$ $\mathbf{indxv}(3)$
\State $\mathbf{e}_i$ $\gets$ 0
\State $\mathbf{e}_j$ $\gets$ 0
\State $\mathbf{e}_k$ $\gets$ 0
\State $\mathbf{e}_i(i)$ $\gets$ 1
\State $\mathbf{e}_j(j)$ $\gets$ 1
\State $\mathbf{e}_k(k)$ $\gets$ 1
\State fr $\gets$  d3mlf($f,\mathbf{e}_i,\mathbf{e}_j,\mathbf{e}_k,
\mathbf{q}$) 
\EndFunction
\end{algorithmic}
\end{algorithm}

\subsubsection{Higher order directional derivatives: case 
$\mathbf{f}:\mathbb{R}^m \to \mathbb{R}^n$.}
From the point of view of a multilinear form, the vector directional 
derivative of a function $\mathbf{f}:\mathbb{R}^m \to \mathbb{R}^n$ is a 
direct generalization of the scalar case, as each component of 
$\mathbf{f}$ is a function $f_{\alpha}:\mathbb{R}^m \to \mathbb{R}$. 
This generalization can be done as follow.

Using the standard basis of $\mathbb{R}^n$ we have
\begin{align}\label{FRmRn}
\mathbf{f} (\mathbf{q})= f_{\alpha} 
(\mathbf{q})\mathbf{e}_{\alpha},~\mathbf{q} \in \mathbb{R}^m~\text{and}~ 
\alpha \in \{1,2,\ldots,n\}.
\end{align}
The 4th order Taylor expansion of this function can be written as
\begin{align}\label{TaylorS4F}\nonumber
&\mathbf{f}(\mathbf{q}+ \mathbf{h}) \;=\; f_\alpha(\mathbf{q})\,
\mathbf{e}_\alpha \;+\; h_i \dpara{f_\alpha(\mathbf{q})}{q_i}\,
\mathbf{e}_\alpha \;+ \nonumber\\
&~~+\frac{h_i h_j}{2!} \dparb{^2 f_\alpha(\mathbf{q})}
{q_i}{q_j}\,\mathbf{e}_\alpha + \frac{h_i h_j h_k}{3!} 
\dparc{^3 f_\alpha(\mathbf{q})}{q_i}{q_j}{q_k} 
\mathbf{e}_\alpha + \nonumber\\
&~~+ \frac{h_i h_j h_k h_l}{4!} 
\dpard{^4 f_\alpha(\mathbf{q})}{q_i}{q_j}{q_k}{q_l}\,\mathbf{e}_\alpha + 
\mathbf{R}(\mathbf{q},\mathbf{h}),
\end{align}
with $h_i$ the $i$-th component of the vector $\mathbf{h}$ and
$\norma{\mathbf{R}(\mathbf{q},\mathbf{h})}$ of the order 
$\norma{\mathbf{h}}^5$. 

The generalization of Eq. (\ref{fqmasev}) for this case is
\begin{align}\label{FvecDMLF4} 
&\mathbf{f}(\mathbf{q}+\epsilon_1\mathbf{v}) =\Big[
\epsilon_0 f_\alpha(\mathbf{q})  + 
\epsilon_1 d_{1\alpha}(\mathbf{v})  + \epsilon_2 d_{2\alpha}
(\mathbf{v},\mathbf{v}) + \nonumber\\
&~~+\epsilon_3 d_{3\alpha}(\mathbf{v},\mathbf{v},
\mathbf{v}) + \epsilon_4 d_{4\alpha}(\mathbf{v},\mathbf{v},\mathbf{v},
\mathbf{v}) \Big] \mathbf{e}_\alpha.
\end{align}

The set of Eqs. (\ref{d1v}--\ref{d4vvvv}) 
becomes now
\begin{align}\label{d1va}
&\hspace*{-0.5cm}\mathbf{d}_{1}(\mathbf{v})=\mathbf{f}(\mathbf{q}+
\epsilon_1\mathbf{v}).\epsilon_1 = v_i \dpara{f_\alpha(\mathbf{q})}
{q_i}\,\mathbf{e}_\alpha\\
\label{d2vva}
&\hspace*{-0.5cm}\mathbf{d}_{2}(\mathbf{v},\mathbf{v})=\mathbf{f}(\mathbf{q}+
\epsilon_1\mathbf{v}).\epsilon_2 = v_i v_j 
\dparb{^2 f_\alpha(\mathbf{q})}{q_i}{q_j}\,\mathbf{e}_\alpha\\
\label{d3vvva}
&\hspace*{-0.5cm}\mathbf{d}_{3}(\mathbf{v},\mathbf{v},\mathbf{v})=\mathbf{f}(\mathbf{q}+
\epsilon_1\mathbf{v}).\epsilon_3 = v_i v_j v_k 
\dparc{^3 f_\alpha(\mathbf{q})} {q_i}{q_j}{q_k}\,\mathbf{e}_\alpha\\
\label{d4vvvva}
&\hspace*{-0.5cm}\mathbf{d}_{4}(\mathbf{v},\mathbf{v},\mathbf{v},\mathbf{v})=
\mathbf{f}(\mathbf{q}+\epsilon_1\mathbf{v}).\epsilon_4 = v_i v_j v_k v_l 
\dpard{^4 f_\alpha(\mathbf{q})}{q_i}{q_j}{q_k}{q_l}\,\mathbf{e}_\alpha.
\end{align}
 As before, the dependence on $\mathbf{f}$ and $\mathbf{q}$ is obviated 
 in the left hand side of these equations. Observe now that the vector 
 space $Y$--see Eq. (\ref{eq:1})--is $\mathbb{R}^n$.

In general, for different arguments we have
\begin{align}\label{d2vec}
&\dparb{^2 f_\alpha(\mathbf{q})}{q_i}{q_j}\,x_i y_j \,\mathbf{e}_\alpha =
\mathbf{d}_2(\mathbf{x},\mathbf{y}) = \frac{1}{2}\Big[
\mathbf{d}_2(\mathbf{x}+\mathbf{y},\mathbf{x}+\mathbf{y}) +\nonumber\\
&~~ -\mathbf{d}_2(\mathbf{x},\mathbf{x}) -\mathbf{d}_2(\mathbf{y},\mathbf{y})\Big]
\end{align}
\begin{align}\label{d3vec}
&\hspace{-0.5cm}\dparc{^3 f_\alpha(\mathbf{q})} {q_i}{q_j}{q_k}\,x_i y_j z_k \,
\mathbf{e}_\alpha=\mathbf{d}_3(\mathbf{x},\mathbf{y},\mathbf{z}) = 
\frac{1}{2}\Big[\mathbf{d}_3(\mathbf{x}+\mathbf{y},\mathbf{x}+\mathbf{y},
\mathbf{z}) +\nonumber\\
&~~- \mathbf{d}_3(\mathbf{x},\mathbf{x},\mathbf{z}) - 
\mathbf{d}_3(\mathbf{y},\mathbf{y},\mathbf{z})\Big],
\end{align}
with 
\begin{align}\label{d3xxzv}
&\hspace{-0.5cm}\mathbf{d}_3(\mathbf{x},\mathbf{x},\mathbf{z})=\frac{1}{6}\Big[
\mathbf{d}_3(\mathbf{z} + \mathbf{x},\mathbf{z} + \mathbf{x},
\mathbf{z} + \mathbf{x}) + \mathbf{d}_3(\mathbf{z} - \mathbf{x},
\mathbf{z} - \mathbf{x},\mathbf{z} - \mathbf{x}) + \nonumber\\
&~~-2 \mathbf{d}_3(\mathbf{z},\mathbf{z},\mathbf{z})\Big],
\end{align}
and for the rank 4 we have
\begin{align} \label{d4vec}
&\hspace{-0.5cm}\dpard{^4 f_\alpha(\mathbf{q})}{q_i}{q_j}{q_k}{q_l}\,x_i y_j z_k w_l
\mathbf{e}_\alpha=\mathbf{d}_4(\mathbf{x},\mathbf{y},\mathbf{z},
\mathbf{w}) = 
\frac{1}{2}\Big[\mathbf{d}_4(\mathbf{x}+\mathbf{y},\nonumber\\
&\mathbf{x}+\mathbf{y},
\mathbf{z},\mathbf{w}) - \mathbf{d}_4(\mathbf{x},\mathbf{x},\mathbf{z},
\mathbf{w})-\mathbf{d}_4(\mathbf{y},\mathbf{y},\mathbf{z},\mathbf{w})\Big],
\end{align}
with
\begin{align}
&\mathbf{d}_4(\mathbf{x},\mathbf{x},\mathbf{z},\mathbf{w})=\frac{1}{2}\Big[
\mathbf{d}_4(\mathbf{x},\mathbf{x},\mathbf{z}+\mathbf{w},\mathbf{z}+
\mathbf{w})-\mathbf{d}_4(\mathbf{x},\mathbf{x},\mathbf{z},\mathbf{z}) +\nonumber\\
&~~-\mathbf{d}_4(\mathbf{x},\mathbf{x},\mathbf{w},\mathbf{w})\Big]
\end{align}
and
\begin{align}
&\mathbf{d}_4(\mathbf{x},\mathbf{x},\mathbf{z},\mathbf{z}) = 
\frac{1}{12}\Big[\mathbf{d_4xz}^{+} + \mathbf{d_4xz}^{-} - 
2 \mathbf{d}_4(\mathbf{x},\mathbf{x},\mathbf{x},\mathbf{x}) +\nonumber\\
&~~-2 \mathbf{d}_4(\mathbf{z},\mathbf{z},\mathbf{z},\mathbf{z})\Big],
\end{align}
with
\begin{eqnarray}
\mathbf{d_4xz}^{+}&=&\mathbf{d}_4(\mathbf{x+z},\mathbf{x+z},
\mathbf{x+z},\mathbf{x+z})\\
\mathbf{d_4xz}^{-}&=&\mathbf{d}_4(\mathbf{x-z},\mathbf{x-z},
\mathbf{x-z},\mathbf{x-z}).
\end{eqnarray}

\begin{algorithm}[htb]
\caption{Third order vector directional derivative}
\label{dmlfr3}
\begin{algorithmic}[1]
\Function{{\bf fr}=d3mlf}{$\mathbf{f},\mathbf{x},\mathbf{y},\mathbf{z},
\mathbf{q},n$}
\State {\bf fr} $\gets$ $\frac{\text{d3qxz}(\mathbf{f},\mathbf{x}+
\mathbf{y},\mathbf{z},\mathbf{q},n) - \text{d3qxz}(\mathbf{f},
\mathbf{x},\mathbf{z},\mathbf{q},n) - 
\text{d3qxz}(\mathbf{f},\mathbf{y},\mathbf{z},\mathbf{q},n)}{2}$
\EndFunction

\Function{{\bf fr}=d3qxz}{$\mathbf{f},\mathbf{x},\mathbf{z},
\mathbf{q},n$}
\State {\bf fr} $\gets$ $\frac{\text{d3q}(\mathbf{f},\mathbf{z}+
\mathbf{x},\mathbf{q},n)+\text{d3q}(\mathbf{f},\mathbf{z}-\mathbf{x},
\mathbf{q},n) - 2\,\text{d3q}(\mathbf{f},\mathbf{z},\mathbf{q},n)}{6}$
\EndFunction

\Function{{\bf fr}=d3q}{$\mathbf{f},\mathbf{v},\mathbf{q},n$}
\State {\bf fqd} $\gets$ $\mathbf{f}(\mathbf{q}+\epsilon_1\,\mathbf{v})$
\State {\bf fr} $\gets$ {\bf fqd}$.\epsilon_3$
\EndFunction
\end{algorithmic}
\end{algorithm}

Formally Eqs. (\ref{d1va}--\ref{d4vec}) are 
the same as their scalar counterpart, and from a computational point of 
view they are essentially equivalent. Nevertheless they allow much more 
concise formulas. For instance, the product of the Jacobian 
matrix of $\mathbf{f(q)}$ with the vector $\mathbf{v}$ is given by 
\begin{align}\label{eqjfv}
[\mathbf{J\,f}(\mathbf{q})]\mathbf{v} = \mathbf{d}_{1}(\mathbf{v}),
\end{align}
which, using Eq. (\ref{d1va}), can be computed with only one evaluation 
of the vector function $\mathbf{f}$ overloaded to accept dual variables. 
Also, using the coordinate representation of this multilinear form we 
obtain the derivatives of all the functions $f_\alpha$ at once. 
Mathematically
\begin{align}\label{dmlfmn1}
\hspace*{-0.5cm}\mathbf{d}_{1}(\mathbf{e}_i) = \dpara{f_{\alpha}(\mathbf{q})}{q_i}\,
\mathbf{e}_\alpha = \bigg[\dpara{f_1(\mathbf{q})}{q_i}~  
\dpara{f_2(\mathbf{q})}{q_i}~ \cdots ~ \dpara{f_n(\mathbf{q})}{q_i}
\bigg]^\text{T}.
\end{align} 
Similarly, for the ranks 2, 3 and 4 we have
\begin{align}
\label{dmlfmn2}
&\hspace*{-0.5cm}\mathbf{d}_{2}(\mathbf{e}_i,\mathbf{e}_j) = \bigg[
\dparb{^2 f_1(\mathbf{q})}{q_i}{q_j}~ 
\dparb{^2 f_2(\mathbf{q})}{q_i}{q_j}~ \cdots ~ 
\dparb{^2 f_n(\mathbf{q})}{q_i}{q_j}
\bigg]^\text{T}\\
\label{dmlfmn3}
&\hspace*{-0.5cm}\mathbf{d}_{3}(\mathbf{e}_i,\mathbf{e}_j,\mathbf{e}_k) = \bigg[
\dparc{^3 f_1(\mathbf{q})}{q_i}{q_j}{q_k}~ 
\dparc{^3 f_2(\mathbf{q})}{q_i}{q_j}{q_k}~ \cdots ~ 
\dparc{^3 f_n(\mathbf{q})}{q_i}{q_j}{q_k}
\bigg]^\text{T}\\
\label{dmlfmn4}
&\hspace*{-0.5cm}\mathbf{d}_{4}(\mathbf{e}_i,\mathbf{e}_j,\mathbf{e}_k,\mathbf{e}_l) = 
\bigg[\dpard{^4 f_1(\mathbf{q})}{q_i}{q_j}{q_k}{q_l}~ 
\dpard{^4 f_2(\mathbf{q})}{q_i}{q_j}{q_k}{q_l}~ \cdots \nonumber\\
&~~\dpard{^4 f_n(\mathbf{q})}{q_i}{q_j}{q_k}{q_l}
\bigg]^\text{T}.
\end{align}

The numerical implementation is practically the same as the scalar case, 
except that $\mathbf{f}$ is now a vector function. For convenience, in 
the Fortran implementation we have included an argument corresponding to 
its dimension $n$. As an example, the version of algorithm 
\ref{dirderr3} for this case is given in algorithm \ref{dmlfr3}.

\section{Numerical examples and comparisons} \label{Ejemplos}
The proposed dual number formulation is coded in Matlab and Fortran
languages. These implementations along the presented examples in this 
section are available at \cite{r_peon2022M}. 

\subsection{Comparison between dual numbers and the FMFD method
to compute a 
directional derivative of order 4}
A function commonly used to test the effectiveness of optimization 
algorithms is the inverted cosine wave function \cite{Penunuri2016}:
\begin{align}
&f(\mathbf{x})=-\sum_{i=1}^{m-1}\exp\left(-\frac{x_i^2+x_{i+1}^2+0.5\,
      x_i\,x_{i+1}}{8}\right) \times \nonumber\\
			&~~\cos
			\left(4\sqrt{x_i^2+x_{i+1}^2+0.5\,x_i\,x_{i+1}}\right).
\label{eq:ICWF}
\end{align}
Table \ref{derdirtab1} shows the computation of the fourth order 
directional derivative for $m\in\{5,7,10,15,20\}$ along the vectors 
$\mathbf{x}=[1~2~\dots~m]^{\text{T}}$, $\mathbf{y}=\sin \mathbf{x}$, 
$\mathbf{z}=\cos \mathbf{x}$, and $\mathbf{w}=\sqrt {\mathbf{x}}$ at the 
point $\mathbf{q}=\log \mathbf{x}$. We consider element-wise 
operation--for instance $\sin([a~b]^{\text{T}})$ $=$ $[\sin a~\sin b]^
{\text{T}}$. The order of approximation used in the FMFD method was 4. 
It is worthwhile to mention that if this order of approximation is 
increased to 8 also the elapsed time considerably increases.

\begin{table}[hbt]
\caption{Comparison between the FMFD method and dual numbers formulation
to compute the four order directional derivative. 
For the FMFD method we used quadruple precision, an order of 
approximation of 4 and a step size $h=10^{-5}$.} 
\centering
\scalebox{0.9}{
\begin{tabular}{l l l l l l}
\toprule
\multirow{2}{*}{$m$} & \multicolumn{2}{c} {Directional derivative}  & &\multicolumn{2}{c} {Time (s)} \\ \cmidrule{2-3} \cmidrule{5-6} 
& FMFD  & Dual numbers & & FMFD & Dual numbers \\ \midrule
5  & 3581.554 & 3581.531 & &1.3    &0.0027 \\ 
7  & 3333.292 & 3333.211 & &7.6   &0.0036 \\ 
10 & 1511.809 & 1511.750 & &47.1   &0.0043 \\ 
15 & 3092.521 & 3092.453 & &375.1  &0.0052 \\ 
20 & 2099.013 & 2098.912 & &1594.0 &0.0064  \\ 
\bottomrule
\end{tabular}
}
\label{derdirtab1}
\end{table}

\begin{figure}[htb]
\hspace*{-0.5cm}\includegraphics[scale=0.95]{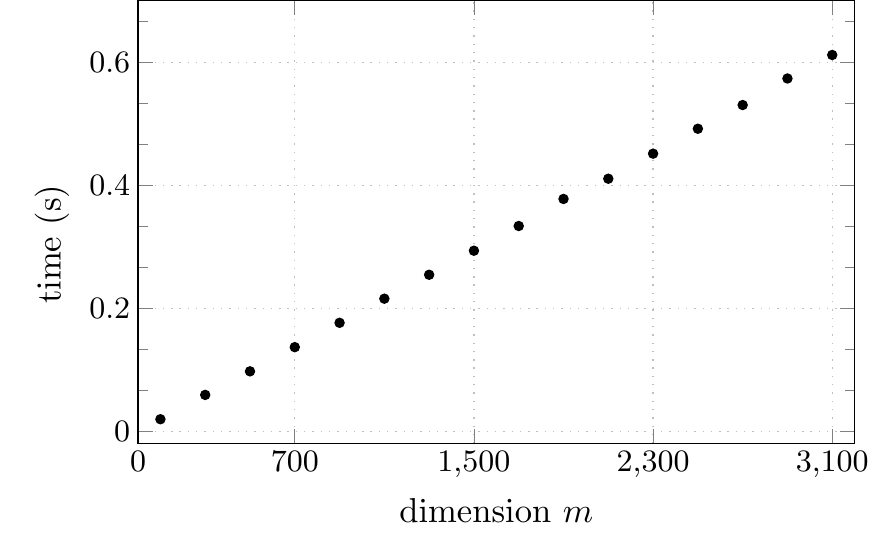} 
\caption{Time complexity for the fourth order directional derivative 
using dual numbers.}
\label{figure1}
\end{figure}

From Table \ref{derdirtab1} we can see how efficient the dual number 
approach is compared to FMFD for this kind of problems. In fact Fig 
\ref{figure1} show a linear time complexity of the 
proposed formulation for computing this directional derivative; the 
complexity using the FMFD method to compute Eq. (\ref{d4mlffor}) for 
this example would be, at the best, of order $O(m^4)$. In that figure the greater value for $m$ 
was 3100. In contrast, the evaluation of a fourth-rank tensor for $m=3100$ would
require computing $9.23521\times 10^{13}$ components, making this approach
prohibitive\footnote{Actually, due to symmetry properties the number of 
independent components is $3855456322025$, a considerable reduction but 
still prohibitive for a symbolic approach or for the FMFD method.}. On the other side with the proposed formulation the elapsed 
time for $m=3100$ was of 0.6 seconds. The computations were made on an 
Intel(R) Core(TM) i7-7700 CPU @ 3.60 GHz, running Windows 10 
and using the gfortran compiler (gcc version 9.2.0). 

\subsection{Calculation of the velocity, acceleration, jerk and 
jounce/snap vectors} \label{secrcre}
Let $\mathbf{r}:\mathbb{R}\to \mathbb{R}^n$ be the position vector of a 
point particle with generalized coordinates, including possibly the 
time, $\mathbf{q}=(q_1(t),q_2(t),\dots,q_m(t))$. Let $f_\alpha:
\mathbb{R}^m \to \mathbb{R}$ be the function which gives the 
$\alpha$-component of $\mathbf{r}$. In the standard basis of 
$\mathbb{R}^n$, the position vector can be written as
\begin{equation}
\mathbf{r}(t) = f_\alpha (\mathbf{q}(t))\,\mathbf{e}_\alpha,
\label{rvec_pos}
\end{equation}

The velocity, acceleration, jerk and jounce vectors of a point particle 
can be obtained differentiating Eq. (\ref{rvec_pos}) as many times as 
required~--if $\mathbf{r}$ depends explicitly on time, set $q_m(t)=t$. 
Such kinematics quantities are given by \cite{Peon2023}
\begin{align}
&\hspace{-0.5cm}\mathbf{\dot{r}} = \dpara{f_\alpha(\mathbf{q})}{q_i} 
\dot{q}_i \, \mathbf{e}_\alpha \label{velocity}\\
&\hspace{-0.5cm}\mathbf{\ddot{r}} = \left[\frac{\partial f_\alpha 
(\mathbf{q})} {\partial q_i}\, \ddot{q}_i +
\frac{\partial^{2} f_\alpha (\mathbf{q})}{\partial q_i \partial q_j}\, 
\dot{q}_i \, \dot{q}_j \right]\mathbf{e}_\alpha \label{acceleration}\\
&\hspace{-0.5cm}\mathbf{\dddot{r}}=\Bigg[\dparc{^3 f_\alpha(\mathbf{q})}{q_i}{q_j}{q_k} 
\dot{q}_i \dot{q}_j \dot{q}_k + 3 \frac{\partial^{2} 
f_\alpha (\mathbf{q})}{\partial q_i \partial q_j}\; \ddot{q}_i \; 
\dot{q}_j \,+ \nonumber\\
&+\frac{\partial f_\alpha(\mathbf{q})}{\partial q_i}\; 
\dddot{q}_i \Bigg]\mathbf{e}_\alpha \label{pulso} \\
&\hspace{-0.5cm}\ddddot{\mathbf{r}} = \Bigg [\dpard{^4 f_\alpha(\mathbf{q})}
{q_i}{q_j}{q_k}{q_l}\dot{q}_i\dot{q}_j \dot{q}_k \dot{q}_l + 
6 \dparc{^3 f_\alpha(\mathbf{q})}{q_i}{q_j}{q_k} \ddot{q}_i \dot{q}_j 
\dot{q}_k +\nonumber\\
&+ 4 \dparb{^2 f_\alpha(\mathbf{q})}{q_i}{q_j} \dddot{q}_i 
\dot{q}_j  \, +  3\dparb{^2 f_\alpha(\mathbf{q})}{q_i}{q_j} 
\ddot{q}_i \ddot{q}_j + \dpara{f_\alpha(\mathbf{q})}{q_i} 
\ddddot{q}_i \Bigg ]\mathbf{e}_\alpha \label{jounce}.
\end{align}
These equations can be written in terms of directional derivatives using 
Eqs. (\ref{d1va}--\ref{d4vec}). They are
\begin{align}
\mathbf{\dot{r}} &= \mathbf{d}_1(\mathbf{\dot{q}}) \label{veld1mlf}\\
\mathbf{\ddot{r}} &= \mathbf{d}_2(\mathbf{\dot{q}},\mathbf{\dot{q}}) +
\mathbf{d}_1(\mathbf{\ddot{q}})  \label{aceld2mlf}\\
\mathbf{\dddot{r}} &= \mathbf{d}_3(\mathbf{\dot{q},\dot{q},\dot{q}}) + 
3 \mathbf{d}_2(\mathbf{\dot{q}},\mathbf{\ddot{q}}) +
\mathbf{d}_1(\mathbf{\dddot{q}}) \label{pulsod3mlf}\\
\ddddot{\mathbf{r}} &= 
\mathbf{d}_4(\mathbf{\dot{q},\dot{q},\dot{q},\dot{q}}) +
6\,\mathbf{d}_3(\mathbf{\dot{q},\dot{q},\ddot{q}}) +
4\, \mathbf{d}_2(\mathbf{\dot{q},\dddot{q}}) + \nonumber\\
&+3\, \mathbf{d}_2(\mathbf{\ddot{q},\ddot{q}}) +
\mathbf{d}_1(\mathbf{\ddddot{q}}) \label{jounced4mlf}.
\end{align}

Observe the conciseness of Eqs. (\ref{veld1mlf}--\ref{jounced4mlf}),  and 
more importantly, these equations are computed accurately and much more 
efficiently than a traditional approach, using Eqs. (\ref{velocity}--\ref{jounce}), where the explicit derivatives and sums are 
required. Indeed, as we have seen in Sec. \ref{sec:DNDD}, all the 
directional derivatives can be computed by \textit{essentially} 
evaluating $\mathbf{f}(\mathbf{q} + \epsilon_1\,\mathbf{v})$ = 
$f_\alpha (\mathbf{q+\epsilon_1\,\mathbf{v}})\,\mathbf{e}_\alpha$.

Another method to compute the above kinematics quantities (at least for 
robots and mechanisms) is screw theory 
\cite{Rico1999, Gallardo2014, Muller2014}, 
which also would give results at the same accuracy as the dual number 
approach. Nevertheless, compared to the dual number formulation, this 
approach is unpractical since each problem would require the construction 
of specific screws. Moreover, nowadays, the screw theory method to 
compute kinematics quantities is devised for kinematic chains, and it is 
not clear how to directly apply it if only the position vector of a point 
particle is known, which contrasts with our formulation. 

Eqs. (\ref{veld1mlf}--\ref{jounced4mlf}) are specially useful when vectors $\mathbf{q}$, 
$\mathbf{\dot{q}}$, $\mathbf{\ddot{q}}$, $\mathbf{\dddot{q}}$, 
$\mathbf{\ddddot{q}}$, are know quantities for a given time $t_0$. In 
the area of robotics and mechanisms this is usually true, but also there 
are cases where instead of knowing these vectors, equations for the 
generalized coordinates (joint variables) as function of time are known. 
This mean that $\mathbf{q} = \mathbf{g}(t)$ is available. In this case 
the dual number approach to compute the kinematics quantities is even 
more simple. Defining the dual quantity
\begin{align}\label{rdt}
\mathbf{\tilde{r}} &= \mathbf{r}(\mathbf{g}(t + \epsilon_1)),
\end{align}
the velocity, acceleration, jerk and jounce/snap are computed by
\begin{align}\label{rdt_1}
\mathbf{\dot{r}}   &= \mathbf{\tilde{r}}.\epsilon_1\\ \label{rdt_2}
\mathbf{\ddot{r}}  &= \mathbf{\tilde{r}}.\epsilon_2\\ \label{rdt_3}
\mathbf{\dddot{r}} &= \mathbf{\tilde{r}}.\epsilon_3\\ \label{rdt_4}
\mathbf{\ddddot{r}}&= \mathbf{\tilde{r}}.\epsilon_4,
\end{align}
respectively. Moreover, Eqs. (\ref{rdt}--\ref{rdt_4})
can also be used to compute all the kinematics 
quantities at a given time $t_0$, even if the function $\mathbf{g}(t)$ 
is unknown, provided $\mathbf{q}$, $\mathbf{\dot{q}}$, 
$\mathbf{\ddot{q}}$, $\mathbf{\dddot{q}}$, $\mathbf{\ddddot{q}}$, are 
know quantities at the time $t_0$. This is achieved using the Taylor 
expansion of the unknown function $\mathbf{g}(t)$ about the point $t_0$.
Furthermore, since what is needed 
to compute Eqs.  (\ref{velocity}--\ref{jounce}) are the quantities $\mathbf{q}$, 
$\mathbf{\dot{q}}$, $\mathbf{\ddot{q}}$, $\mathbf{\dddot{q}}$, 
$\mathbf{\ddddot{q}}$ and not the explicit time dependence of 
$\mathbf{q} = \mathbf{g}(t)$, the obtained values of 
such kinematics quantities at the time $t_0$ are exact.

\subsubsection{A hypothetical example}
An interesting AD property is that a closed form expression for the 
function to derive is not necessary, but rather an algorithmic way to 
compute it. For example, 
let $\mathbf{f}:\mathbb{R}^2 \to \mathbb{R}^3$ given by
\begin{equation}
\mathbf{f}(q_1,q_2) = [\log(q_1\,q_2^2)~~~\sqrt{q_2}/q_1~~~ 
\underbrace{\sin(\sin( \dots \sin(q_1\,q_2)\dots)}_{\mbox{$100$ 
times}}]^{\text{T}},
\label{eq:r12ex}
\end{equation}
this function needs to nest  the sine function 100 times, which is 
unpractical. Instead we can write a do-loop to code this function.

Assuming Eq. (\ref{eq:r12ex}) is the position vector of a point particle
and that $q_1$ and $q_2$ are generalized coordinates depending on time, 
Table \ref{compE1} shows the velocity, acceleration, jerk and jounce for
$q_1=1.1$, $q_2=2.2$, 
$\mathbf{\dot{q}} = [0.5~~-2.7]^{\text{T}}$, 
$\mathbf{\ddot{q}} = [-0.1~~0.7]^{\text{T}}$,
$\mathbf{\dddot{q}} = [0.3~~0.5]^{\text{T}}$, and
$\mathbf{\ddddot{q}} = [-0.2~~0.1]^{\text{T}}$.
Notice that a screw theory approach for this case is not possible or at 
least it cannot be applied directly. A symbolic approach is highly 
inefficient while a finite differences approach is, in general, both; 
inaccurate and inefficient.
\begin{table}[htb]
\caption{Kinematic quantities computed with dual numbers}
\centering
\begin{tabular}{lr}
\toprule
Kinematic Quantity                     & Components values  \\
\midrule 								
\multirow{3}{*}{}                      & $-2.0000$ \\
Velocity ~\normalsize{[L]/[t]}         & $-1.4403$ \\
                                       & $ 0.0207$ \\
\midrule													
\multirow{3}{*}{}                      & $-2.6736$ \\
Acceleration ~\normalsize{[L]/[t]$^2$} & $ 1.1388$ \\
                                       & $-0.1441$ \\
\midrule													
\multirow{3}{*}{}                      & $-4.0120$ \\
Jerk ~\normalsize{[L]/[t]$^3$}         & $-2.7000$ \\
                                       & $ 1.0395$ \\
\midrule													
\multirow{3}{*}{}                      & $-15.1909$ \\
Jounce/Snap ~\normalsize{[L]/[t]$^4$}  & $ 6.2371$ \\
                                       & $-8.0518$ \\
\bottomrule													
\end{tabular} \label{compE1}
\end{table}

\subsubsection{A practical example}
\begin{figure}[htb]
\begin{center}
\includegraphics[scale=0.25]{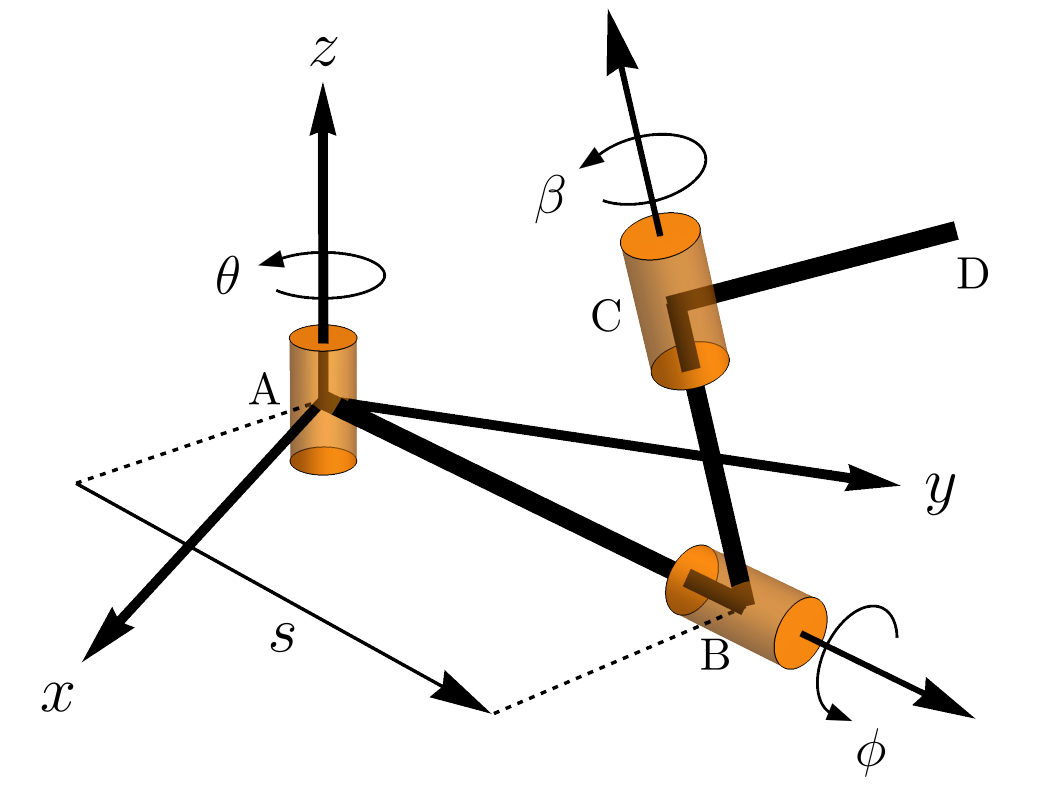} 
\caption{The RCR robot manipulator.}
\label{figure2}
\end{center}
\end{figure}

Let us consider the RCR robot manipulator of Fig. \ref{figure2}, with A, 
B, C the rotational, cylindrical and rotational joints respectively, with 
$\overline{BC}=3$ and $\overline{CD}=2$. 
We are interested in calculating the velocity, acceleration, jerk and 
jounce/snap for a given time where 
\begin{align}
\mathbf{q} &= [\theta=\pi/2~~ \phi=0~~ s=2~~ \beta=0]^{\text{T}},\nonumber \\
\mathbf{\dot{q}} &= [1~~ 5~~ 1~~ 1]^{\text{T}},\nonumber\\
\mathbf{\ddot{q}} &= [1~~ 0~~ 2~~ 1]^{\text{T}},\nonumber\\
\mathbf{\dddot{q}} &= [1~~ 2~~ 3~~ 4]^{\text{T}},\nonumber\\
\mathbf{\ddddot{q}} &= [4~~ 5~~ 6~~ 7]^{\text{T}}.
\label{eqsqsvals}
\end{align}

The position vector of the end effector, point D, can be found using the 
homogeneous transformation matrix method \cite{Stejskal1996}. This 
position vector is obtained taking the ﬁrst three elements of the 4-th 
column of the matrix $\mathbf{R_{\text{D}}}$  given by
\begin{equation}
\mathbf{R_{\text{D}}} = \mathbf{T}_{1}\, \mathbf{T}_{2}\, 
\mathbf{T}_{3}\, \mathbf{T}_{4}\, \mathbf{T}_5,
\end{equation}
where
\begin{align}
\mathbf{T}_1 &= \mathbf{HTM}(\theta,[0,0,1],[0,0,0]) \nonumber \\
\mathbf{T}_2 &= \mathbf{HTM}(\phi,[1,0,0],[s,0,0])\nonumber \\
\mathbf{T}_3 &= \mathbf{HTM}(\beta,[0,0,1],[0,0,0])\nonumber \\
\mathbf{T}_4 &= \mathbf{HTM}(0,[0,0,1],[0,0,BC])  \nonumber \\
\mathbf{T}_5 &= \mathbf{HTM}(0,[1,0,0],[CD,0,0]),
\label{T12345}
\end{align}
being $\mathbf{HTM}(\psi,\mathbf{\hat{u}},\mathbf{S})$ the homogeneous 
transformation matrix of angle $\psi$ about the axis $\mathbf{\hat{u}}$
and for the translation vector $\mathbf{S}$. From Eqs (\ref{veld1mlf}--\ref{jounced4mlf}) 
or Eqs. (\ref{rdt_1}--\ref{rdt_4}), the mentioned 
kinematic quantities can be efficiently computed once the involved 
operators and functions are overloaded to accept dual variables as 
arguments. The numerical values are shown in Table \ref{compE2}. 
It is worthwhile to mention that these values coincide with those 
computed using screw theory (\ref{A1}).
\begin{table}[htb]
\caption{Kinematic quantities for the RCR robot manipulator, 
computed with dual numbers.}
\centering
\begin{tabular}{lr}
\toprule
 Kinematic Quantity                           & Component values  \\
\midrule
\rowcolor[HTML]{ECF4FF} 
\cellcolor[HTML]{ECF4FF}                      & $9.0000$ \\
\rowcolor[HTML]{ECF4FF} 
\cellcolor[HTML]{ECF4FF} 
Velocity ~\normalsize{[L]/[t]}                & $1.0000$  \\
\rowcolor[HTML]{ECF4FF} 
\multirow{-3}{*}{\cellcolor[HTML]{ECF4FF}}    & $ 0.0000$ \\
                                              & $-8.0000$ \\
Acceleration ~\normalsize{[L]/[t]$^2$}        & $22.0000$ \\
\multirow{-3}{*}{}                            & $-55.0000$ \\
\rowcolor[HTML]{ECF4FF} 
\cellcolor[HTML]{ECF4FF}                      & $-267.0000$ \\
\rowcolor[HTML]{ECF4FF} 
\cellcolor[HTML]{ECF4FF}   
Jerk ~\normalsize{[L]/[t]$^3$}                & $15.0000$ \\
\rowcolor[HTML]{ECF4FF} 
\multirow{-3}{*}{\cellcolor[HTML]{ECF4FF}}    & $30.0000$ \\
                                              & $189.0000$\\
Jounce/Snap ~\normalsize{[L]/[t]$^4$}         & $-978.0000$\\
\multirow{-3}{*}{}                            & $891.0000$\\
\bottomrule
\end{tabular}\label{compE2}
\end{table}


\subsection{Examples of partial derivatives}
It is illustrative to compute some partial derivatives using the 
proposed approach as well as using the FMFD method. This approach 
however, is not recommended to compute general directional derivatives.
As an example, Table \ref{EMPD} shows some partial derivatives
of the function 
\begin{align}
f(x,y,z,w,u) &= \cos(x\,y/u)\,z/w + 3\,\sin(x\,u)\,\sin(y/u) \times \nonumber \\
&\log(x\,y\,z/(u\,w)),
\label{eq:eff}
\end{align}
evaluated at the point $\mathbf{x}_0=(1.1,2.2,3.3,4.4,5.5)$ using the 
FMFD method and dual numbers.
\newsavebox\DPA
\begin{lrbox}{\DPA}
  \begin{minipage}[b]{2.2cm}
    \begin{equation*}
		\dpara{f(\mathbf{x}_0)}{x}
    \end{equation*} 
  \end{minipage}
\end{lrbox}

\newsavebox\DPB
\begin{lrbox}{\DPB}
  \begin{minipage}[b]{2.2cm}
    \begin{equation*}
		\dparb{^2 f(\mathbf{x}_0)}{z\,}{z}
    \end{equation*} 
  \end{minipage}
\end{lrbox}

\newsavebox\DPC
\begin{lrbox}{\DPC}
  \begin{minipage}[b]{2.2cm}
    \begin{equation*}
		\dparc{^3 f(\mathbf{x}_0)}{u\,}{w\,}{y}
    \end{equation*} 
  \end{minipage}
\end{lrbox}

\newsavebox\DPD
\begin{lrbox}{\DPD}
  \begin{minipage}[b]{2.2cm}
    \begin{equation*}
		\dpard{^4 f(\mathbf{x}_0)}{u\,}{z\,}{w\,}{x\,}
    \end{equation*} 
  \end{minipage}
\end{lrbox}

\newsavebox\DNFR
\begin{lrbox}{\DNFR}
  \begin{minipage}[b]{2.5cm}
    \small{Dual number\\
		formulation}
  \end{minipage}
\end{lrbox}

\begin{table*}[h]
\caption{Examples of partial derivatives.}
\centering
\scalebox{1.0}{
\begin{tabular}{l l l l l}
\toprule
\multirow{2}{*}{Partial derivative} & \multirow{2}{*}{\scalebox{0.9}
{\usebox{\DNFR}}}
& \multicolumn{3}{c}  {Finite differences to 8th order 
$\left(h=10^{-5}\right)$} \\ \cmidrule{3-5} 
&   & single precision  & double precision & quadruple precision \\ \midrule					
\usebox{\DPA} & $-7.3040034203$ & $-5.64597130$ & $-7.3040034219$ &$-7.3040034203$ \\ 
\usebox{\DPB} & ~~\,$0.0247895122$ & $-40885.4766$ & ~~\,$0.0248116774$  &~~\,$0.0247895122$  \\ 
\usebox{\DPC} & $-0.1312934721$ & ~~\,$198403296$   & -2.4530159385   &$-0.1312934721$ \\ 
\usebox{\DPD} & $-0.0030955304$ & $-3.06289972\times10^{13}$ & $64117.282468$ & 
$-0.0030955304$ \\ 
& & & & \\
\bottomrule
\end{tabular}
}
\label{EMPD}
\end{table*}

This table shows that the FMFD method fails to compute the derivatives 
if single and double precision variables are used. Nevertheless, using 
quadruple precision, the results coincide with the dual number 
formulation which was coded to double precision.

\section{Conclusions}\label{conclusions}
Dual numbers allow a precise and efficient calculation of higher order 
directional derivatives. The evaluation of a scalar function of vector 
variable at the dual point $\mathbf{q} + \epsilon_1\,\mathbf{v}$ permits 
the computation of all the directional derivatives along a single vector 
$\mathbf{v}$, at the point $\mathbf{q}$. Using this result, formulas for
directional derivatives along different vectors are obtained. 
Directional derivatives are special cases of multilinear forms; thus, 
with the coordinate representation of these multilinear forms, mixed 
partial derivatives can be computed. The efficiency of the 
proposed formulation is remarkable. In one of the presented examples, 
the fourth order directional derivative of a function with 3100 
variables is obtained in less than a second. A brute force computation 
of this directional derivative  using finite differences would require 
the computation of a tensor of $9.23521 \times 10^{13}$ components. 
However, our numerical experiments show that finite differences can be 
an excellent option to compute partial derivatives when the involved
dimensions are moderated. Although, it may require variables of high 
precision (quadruple precision or above). For a vector function of 
several variables, the evaluation in the aforementioned dual point 
permits the computation of all the vector directional derivatives 
(directional derivatives for a function $\mathbf{f}:\mathbb{R}^m\to 
\mathbb{R}^n$). With these vector directional derivatives, formulas to 
efficiently compute kinematic quantities can be deduced.  When the 
generalized coordinates are given as a function $\mathbf{g}(t)$, the 
dual number approach is even more direct since all the kinematic 
quantities can be computed by just evaluating the position vector on 
this function. If this function is unknown, and the set  $\{\mathbf{q,
~\dot{q},~\ddot{q},~\dddot{q},~\ddddot{q}}\}$ is given instead, a Taylor 
series of $\mathbf{g}(t)$ can be used. From a theoretical point of view, 
it would be interesting to construct the rules governing the time 
differentiation of the involved directional derivatives on computing 
kinematics quantities. A close inspection to the arguments of equations 
(\ref{veld1mlf}--\ref{jounced4mlf}) --or in general, Eqs. 
(\ref{d1va}--\ref{d4vec})-- suggests that such rules, and 
generalizations to higher orders, could be found following a combinatorics 
approach. We hope any researcher could shed some light on this respect 
in future studies.

\section*{Acknowledgements}
The authors acknowledge the support of the Consejo Nacional de
Ciencia y Tecnolog\'ia (National Council of Science and Technology,
CONACYT), of M\'exico, through SNI (National Network of Researchers)
fellowships and scholarships.

\section*{CRediT authorship contribution statement}
\textbf{R. Pe\'on-Escalante:} Investigation, Visualization, Writing - 
original draft. 
\textbf{K.B. Cant\'un-Avila:} Conceptualization, Methodology, Software, 
Validation, Writing - original draft, Writing - review 
\& editing draft. 
\textbf{O. Carvente:} Investigation, Validation, Visualization. 
\textbf{A. Espinosa-Romero:} Methodology, Software, 
Validation,  Writing - original draft. 
\textbf{F. Pe\~nu\~nuri:} Conceptualization, Formal Analysis, Methodology, 
Software, Supervision, Validation, Writing - original draft, Writing - 
review \& editing draft.

\section*{Declaration of competing interest}
The authors declare that they have no known competing financial 
interests or personal relationships that could have appeared to
influence the work reported in this paper.

\appendix
\section{Computation of the kinematic quantities for the RCR robot 
manipulator using infinitesimal screw theory}
\label{A1}
This appendix presents the essential elements 
of screw theory to compute the kinematics quantities of the example 
 in sec. \ref{secrcre}. The interested reader can see ref. 
\cite{Custodio2017} for details. We recall that in this section, the 
Einstein's convention about sum over repeated index is not used, and all 
the involved sums are written explicitly.

\begin{figure}[h]
\begin{center}
\includegraphics[scale=0.2]{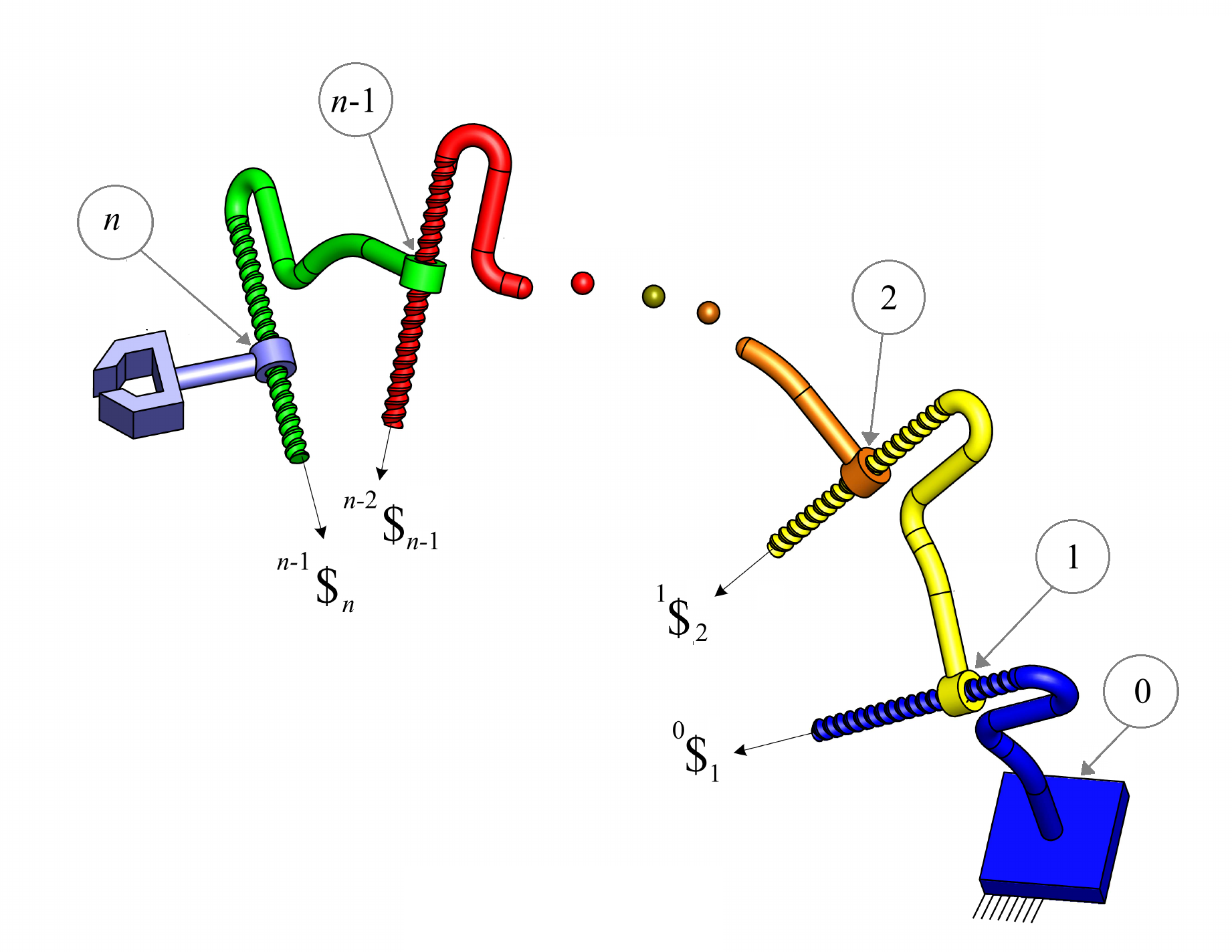} 
\caption{Serial kinematic chain
\cite{Custodio2017}.}
\label{CadenaCinematica}
\end{center}
\end{figure}

Considering the kinematic chain in Fig. \ref{CadenaCinematica}, we will 
say that we have $n + 1$ \emph{bodies} numbered from $0$ to $n$, and the body 
tagged with 0 is the frame of a kinematic chain. We assume that a body $i$ has a
frame attached to it whose origin measured in the frame (or body) $j$ is 
located by the vector $^j\mathbf{r}_i$; from a practical point of view, 
this position vector would locate the joints. From now on, the symbol 
$^a\star_b$ will means that the quantity $\star$ is related to frame $b$ 
measured in frame $a$. For instance, $^j\mathbf{\dot{r}}_n=\,^j
\mathbf{v}_n$ would be the velocity of the origin of the frame $n$ (or just the
velocity of frame $n$) with respect to the frame $j$.

Let us consider a point $p$ on frame (or body) $n$, then its position 
vector on a frame $j$ is ${}^j\mathbf{r}^p_{n}$. Notice this position 
vector is not necessarily equal to ${}^j\mathbf{r}_{n}$, except if the point $p$
is in the origin of frame $n$. To obtain the 
kinematics quantities of this point $p$ using screw theory, the 
so-called \emph{reduced velocity state}, \emph{reduced acceleration state}, 
\emph{reduced jerk state} and  \emph{reduced jounce state} are required. 
These \emph{kinematic states} are defined as:
\begin{align}
{^j\mathbf{V}_n} &= \left[ \begin{array}{c} {^j\gbf\omega_n} \\ 
{^j\mathbf{v}^*_n} \end{array} \right] \label{EQVS}\\
{^j\mathbf{A}_n} &= \left[ \begin{array}{c} {^j\gbf\alpha_n} \\ 
{^j\mathbf{a}^*_n} \end{array} \right]\label{EQAS}\\
{^j\mathbf{J}_n} &= \left[ \begin{array}{c} {^j\gbf\rho_n} \\ 
{^j\mathbf{j}^*_n} \end{array} \right]\label{EQJKS}\\
{^j\mathbf{S}_n} &= \left[ \begin{array}{c} {^j\gbf\sigma_n} \\ 
{^j\mathbf{s}^*_n} \end{array} \right]\label{EQJSS}.
\end{align}

From the algebraic point of view, these expressions are $6\times1$ 
matrices, although, from a programming point of view, arrays of 6 
elements can be used. The first vector components of the above kinematic 
states correspond to angular velocity, angular acceleration, angular 
jerk, and angular jounce/snap, respectively, of the frame $n$, with 
respect to frame $j$. The second components are related\footnote{Notice 
that the velocity of frame $n$ with respect to frame $j$ is not the 
\emph{reduced velocity} $^{j}\mathbf{v}^*_{n}$ but
$^{j}\mathbf{v}_{n} = \;^{j}\mathbf{v}^*_{n} +\, {}^j\gbf\omega_n \times
{}^j\mathbf{r}_{n}\label{velFn}.$} to their linear 
counterparts, namely, velocity, acceleration, jerk, and jounce/snap. 
Such relations, for the aforementioned point $p$, are:
\begin{align}
^{j}\mathbf{v}^p_{n} &= \;^{j}\mathbf{v}^*_{n} +\, {}^j\gbf\omega_n \times
{}^j\mathbf{r}^p_{n}\label{velSCR}\\
^{j}\mathbf{a}^p_{n} &= \;^{j}\mathbf{a}^*_{n} +\, {}^j\gbf\alpha_n \times
{}^j\mathbf{r}^p_{n} +\, {}^j\gbf\omega_n \times {}^j\mathbf{v}^p_{n} 
\label{accelSCR}\\
^{j}\mathbf{j}^p_{n} &= \;^{j}\mathbf{j}^*_{n} +\, {}^j\gbf\rho_n \times
{}^j\mathbf{r}^p_{n} +\, 2\,{}^j\gbf\alpha_n \times {}^j\mathbf{v}^p_{n}+
{}^j\gbf\omega_n \times {}^j\mathbf{a}^p_{n}\label{jerkvelSCR}\\
^{j}\mathbf{s}^p_{n} &= \;^{j}\mathbf{s}^*_{n} +\, {}^j\gbf\tau_n \times
{}^j\mathbf{r}^p_{n} +\, 3\,{}^s\gbf\rho_n \times {}^j\mathbf{v}^p_{n}+
3\,{}^j\gbf\alpha_n \times {}^j\mathbf{a}^p_{n} +\nonumber\\
&+\,{}^j\gbf\omega_n \times {}^j\mathbf{j}^p_{n}. \label{jouncesnapSCR}
\end{align}

To calculate Eqs. (\ref{EQVS}--\ref{EQJSS}) for the RCR robot 
manipulator (see Fig. \ref{figure2}), we need the associated screws for 
each joint. The most general joint of this robot is the cylindrical 
joint; whose associated screw is the sum of the screws associated to 
rotational and prismatic (translational) joints. Nevertheless, to keep 
the discussion about how to compute the kinematics states as clear as 
possible, we define the \emph{vector} screw 
\begin{align}
{}^j\gbf{\$}_{i}= [^j\$^{\text{R}}_{i}~~^j\$^{\text{T}}_{i}],
\end{align}
whose components are the screws associated to the rotational and 
prismatic joints; given by
\begin{align}
{}^j\$^{\text{R}}_{i}&=\begin{bmatrix}
{}^j\mathbf{e}_{i} \\
{}^j\mathbf{r}_{i} \, \times \, ^j\mathbf{e}_{i}
\end{bmatrix} \text{ and}\\
{}^j\$^{\text{T}}_{i}&=\begin{bmatrix}
\mathbf{0} \\
{}^j\mathbf{e}_{i}
\end{bmatrix},
\end{align}
respectively. As before, $^j\mathbf{r}_i$ locates the position of the 
joints (the origin of a frame) and $^j\hat{\mathbf{e}}_i$ a unit vector 
in the same direction of axis of joint $i$ with respect to frame $j$. 
We also define the following vectors storing the joint variables 
(generalized coordinates)
\begin{align}
\mathbf{q}_{\text{R}}&=[\theta~~0]\\
\mathbf{q}_{\text{T}}&=[0~~s]\\
\mathbf{q}_{\text{C}}&=[\theta~~s].
\end{align}
Since, for instance, products of the form
\begin{align}
\mathbf{\dot{q}}_{\text{R}}\, \cdot\, {}^j\gbf{\$}_{i} &= 
\dot{\theta}\,{}^j\$^{\text{R}}_{i} \\
\mathbf{\dot{q}}_{\text{T}}\, \cdot\, {}^j\gbf{\$}_{i} &= 
\dot{s}\,{}^j\$^{\text{T}}_{i}\\
\mathbf{\dot{q}}_{\text{C}}\, \cdot\, {}^j\gbf{\$}_{i} &= 
\dot{\theta}\,{}^j\$^{\text{R}}_{i} + \dot{s}\,{}^j\$^{\text{T}}_{i}
\end{align}
will appear when computing the kinematic 
states, we define the matrix (or if prefer, vector of vectors)
\begin{align}
 \mathbf{q}&=[\mathbf{q}_{\text{R}}~~\mathbf{q}_{\text{T}}~~\mathbf{q}_{\text{C}}],
\end{align}
and its time derivatives, for example $d\mathbf{q}/dt = 
\mathbf{\dot{q}}=[\mathbf{\dot{q}}_{\text{R}}~~
\mathbf{\dot{q}}_{\text{T}}~~\mathbf{\dot{q}}_{\text{C}}]$. 
Since we are interested in computing the kinematic quantities in the 
\emph{laboratory} system we will use $j=0$. The velocity state for
point $p$ using infinitesimal screws is given by
\begin{align}
 \left[ \begin{array}{c} {^0\gbf\omega_n} \\ 
{^0\mathbf{v}^*_n} \end{array} \right] = \sum_{k=1}^n \mathbf{\dot{q}}_k\,
\cdot\,{}^0\gbf{\$}_{k}.
\end{align}
Bear in mind that the index $k$ in $\mathbf{\dot{q}}_k$ means the $k$-th 
component of $\mathbf{\dot{q}}$ (which is a vector), while in 
${}^0\gbf{\$}_{k}$ is related to frame (or body) $k$.

To calculate the acceleration state we recall the Lie product between two
screws 
\begin{align}
 \$_1=\left[ \begin{array}{c} 
 s_1 \\ 
s_{01} \end{array} \right], \text{ and } 
\$_2=\left[ \begin{array}{c} 
 s_2 \\ 
s_{02} \end{array} \right]
\end{align}
as 
\begin{align}
 [\$_1,\$_2] = \left[ \begin{array}{c} 
 s_1 \times s_2 \\ 
s_1 \times s_{02} - s_2 \times s_{01} \end{array} \right].
\end{align}
and the Lie Screw ${}^0\$\text{L}_n$ defined as
\begin{align}
 {}^0\$\text{L}_n=\sum_{i=1}^{n-1}\sum_{k=i+1}^{n}\,
 [\mathbf{\dot{q}}_i \cdot {}^0\gbf{\$}_{i} ,
  \mathbf{\dot{q}}_k \cdot {}^0\gbf{\$}_{k}],
\end{align}
with this, the acceleration state is given by
\begin{align}
 \left[ \begin{array}{c} {^0\gbf\alpha_n} \\ 
{^0\mathbf{a}^{*}_n} \end{array} \right] = \sum_{k=1}^n \mathbf{\ddot{q}}_k\,
\cdot\,{}^0\gbf{\$}_{k} + {}^0\$\text{L}_n.
\end{align}

To prevent the possibility of typos it is convenient to define the 
functions
\begin{align}
\hspace{-0.5cm} \text{L}02(\mathbf{x},\mathbf{y},a,b)&=
[\mathbf{x}_a \cdot {}^0\gbf{\$}_{a} ,
 \mathbf{y}_b \cdot {}^0\gbf{\$}_{b}] \label{eqSL02}\\
\hspace{-0.5cm}\text{L}03(\mathbf{x},\mathbf{y},\mathbf{z},a,b,c)&=
	[\mathbf{x}_a \cdot {}^0\gbf{\$}_{a} ,
	\text{L}02(\mathbf{y},\mathbf{z},b,c)]\label{eqSL03}\\
\hspace{-0.5cm}\text{L}04(\mathbf{x},\mathbf{y},\mathbf{z},\mathbf{w},a,b,c,d)&=
	[\mathbf{x}_a \cdot {}^0\gbf{\$}_{a} ,
	\text{L}03(\mathbf{y},\mathbf{z},\mathbf{w},b,c,d)] \label{eqSL04},
\end{align}
with this, the Lie screw can be rewritten as
\begin{align}
 {}^0\$\text{L}_n=\sum_{i=1}^{n-1}\sum_{k=i+1}^{n}\,
 \text{L}02(\mathbf{\dot{q}},\mathbf{\dot{q}},i,k).
\end{align}

The jerk state is 
\begin{align}
\left[ \begin{array}{c} {^j\gbf\rho_n} \\ 
{^j\mathbf{j}^*_n} \end{array} \right]=\sum_{k=1}^n \mathbf{\dddot{q}}_k\,
\cdot\,{}^0\gbf{\$}_{k} + {}^0\$\text{J}_n,
\end{align}
with ${}^0\$\text{J}_n$ the jerk screw given by
\begin{align}
 &{}^0\$\text{J}_n =\sum_{i=1}^{n-1}\sum_{k=i+1}^{n}\,
 2\,\text{L}02(\mathbf{\dot{q}},\mathbf{\ddot{q}},i,k) \;+\nonumber\\
&~+ \text{L}02(\mathbf{\ddot{q}},\mathbf{\dot{q}},i,k)\,+\,
 \text{L}03(\mathbf{\dot{q}},\mathbf{\dot{q}},\mathbf{\dot{q}},i,i,k) ,+\,\nonumber\\
 &~+\,\sum_{l=1}^{n-2}\;\sum_{i=l+1}^{n-1}\;\sum_{k=i+1}^{n} 
  2\,\text{L}03(\mathbf{\dot{q}},\mathbf{\dot{q}},\mathbf{\dot{q}},l,i,k).
\end{align}

The jounce/snap state is
\begin{align}
\left[ \begin{array}{c} {^j\gbf\sigma_n} \\ 
{^j\mathbf{s}^*_n} \end{array} \right]=\sum_{k=1}^n \mathbf{\ddddot{q}}_k\,
\cdot\,{}^0\gbf{\$}_{k} + {}^0\$\text{S}_n,
\end{align}
where ${}^0\$\text{S}_n$ is the jounce/snap screw given by
\begin{align}
{}^0\$\text{S}_n&= \text{SIK}+\text{SLIK}+\text{SHLIK}
\end{align}
with
\begin{align}
&\hspace{-0.5cm}\text{SIK}_{}=\sum_{i=1}^{n-1}\sum_{k=i+1}^{n}\,\Big[
3\,\text{L}02(\mathbf{\dot{q}},\mathbf{\dddot{q}},i,k) \,+\,
3\,\text{L}02(\mathbf{\ddot{q}},\mathbf{\ddot{q}},i,k)\,+\nonumber\\
&\hspace{-.0cm}+\text{L}02(\mathbf{\dddot{q}},\mathbf{\dot{q}},i,k)\,+\,   
3\,\text{L}03(\mathbf{\dot{q}},\mathbf{\dot{q}},\mathbf{\ddot{q}},i,i,k)\,+\, \nonumber\\
&+\,2\,\text{L}03(\mathbf{\ddot{q}},\mathbf{\dot{q}},\mathbf{\dot{q}},i,i,k)\,+\,
\text{L}03(\mathbf{\dot{q}},\mathbf{\ddot{q}},\mathbf{\dot{q}},i,i,k)\,+\nonumber\\
&+\text{L}04(\mathbf{\dot{q}},\mathbf{\dot{q}},\mathbf{\dot{q}},\mathbf{\dot{q}},i,i,i,k)
\Big]
\end{align}
\begin{align}
&\hspace{-0.5cm}\text{SLIK}_{} =\sum_{l=1}^{n-2}\sum_{i=l+1}^{n-1}\sum_{k=i+1}^{n}\,\Big[
6\,\text{L}03(\mathbf{\dot{q}},\mathbf{\dot{q}},\mathbf{\ddot{q}},l,i,k)\,+\nonumber\\
&+3\,\text{L}03(\mathbf{\dot{q}},\mathbf{\ddot{q}},\mathbf{\dot{q}},l,i,k)\,+\,
3\,\text{L}03(\mathbf{\ddot{q}},\mathbf{\dot{q}},\mathbf{\dot{q}},l,i,k)\,+\, \nonumber \\
&+\,3\,\text{L}04(\mathbf{\dot{q}},\mathbf{\dot{q}},\mathbf{\dot{q}},\mathbf{\dot{q}},l,l,i,k)\,+\,
3\,\text{L}04(\mathbf{\dot{q}},\mathbf{\dot{q}},\mathbf{\dot{q}},\mathbf{\dot{q}},l,i,i,k)
\Big] 
\end{align}
\begin{align}
\hspace{-0.5cm}\text{SHLIK}_{}&=\sum_{h=1}^{n-3}\;\sum_{l=h+1}^{n-2}\;\sum_{i=l+1}^{n-1}\;\sum_{k=i+1}^{n}\,
3\,\text{L}04(\mathbf{\dot{q}},\mathbf{\dot{q}},\mathbf{\dot{q}},\mathbf{\dot{q}},h,l,i,k) 
\end{align}

Since we can compute all the kinematics states, setting $j=0$ and $n=3$ 
in Eqs. (\ref{velSCR}--\ref{jouncesnapSCR}) and taking $\mathbf{\dot{q}},\mathbf{\ddot{q}},
\mathbf{\dddot{q}},\mathbf{\ddddot{q}}$ as given in 
Eqs. (\ref{eqsqsvals}) with  $\mathbf{T}_k$ $k \in \{1,...5\}$ as in
Eqs. (\ref{T12345}), we obtain the kinematics quantities of table 
\ref{compE2} when
\begin{align}
{}^0\mathbf{e}_{1}&=[0~~0~~1]^{\text{T}}\\
{}^0\mathbf{e}_{2}&=\dfrac{{}^0\mathbf{r}_{2}}{\norma{{}^0\mathbf{r}_{2}}}\\
{}^0\mathbf{e}_{3}&=\dfrac{{}^0\mathbf{r}_{3} - {}^0\mathbf{r}_{2}}
{\norma{{}^0\mathbf{r}_{3} - {}^0\mathbf{r}_{2}}}
\end{align}
with
\begin{align}
{}^0\mathbf{r}_{1}&=[0~~0~~0]^{\text{T}}\\
{}^0\mathbf{r}_{2}&=[\mathbf{T}_{1}\, \mathbf{T}_{2}] (1:3,4)\\
 {}^0\mathbf{r}_{3}&=[\mathbf{T}_{1}\, \mathbf{T}_{2}\, 
\mathbf{T}_{3}\, \mathbf{T}_{4}] (1:3,4),
\end{align}
where $\mathbf{A}(1:3,4)$ means the first three elements of the fourth 
column of matrix $\mathbf{A}$.

\phantomsection
\bibliographystyle{elsarticle-num}


\begin{thebibliography}{10}
\expandafter\ifx\csname url\endcsname\relax
  \def\url#1{\texttt{#1}}\fi
\expandafter\ifx\csname urlprefix\endcsname\relax\def\urlprefix{URL }\fi
\expandafter\ifx\csname href\endcsname\relax
  \def\href#1#2{#2} \def\path#1{#1}\fi

\bibitem{Margossian2019}
C.~C. Margossian, A review of automatic differentiation and its efficient
  implementation, WIREs Data Mining and Knowledge Discovery 22~(4) (2019).
\newblock \href {https://doi.org/10.1002/widm.1305}
  {\path{doi:10.1002/widm.1305}}.

\bibitem{Atilim2018}
A.~G. Baydin, B.~A. Pearlmutter, A.~A. Radul, J.~M. Siskind,
  \href{http://jmlr.org/papers/v18/17-468.html}{Automatic differentiation in
  machine learning: a survey}, Journal of Machine Learning Research 18~(153)
  (2018) 1--43.
\newline\urlprefix\url{http://jmlr.org/papers/v18/17-468.html}

\bibitem{Penunuri2020}
F.~{Pe\~nu\~nuri}, R.~{Pe\'on}-Escalante, D.~{Gonz\'alez}-{S\'anchez},
  M.~{Escalante Soberanis}, Dual numbers and automatic differentiation to
  efficiently compute velocities and accelerations, Acta Applicandae
  Mathematicae 170 (2020) 649--659.
\newblock \href {https://doi.org/10.1007/s10440-020-00351-9}
  {\path{doi:10.1007/s10440-020-00351-9}}.

\bibitem{Chandrasekhar2021}
C.~Aaditya, S.~Saketh, S.~Krishnan, {AuTO}: a framework for automatic
  differentiation in topology optimization, Structural and Multidisciplinary
  Optimization 64 (2021) 4355--4365.
\newblock \href {https://doi.org/10.1007/s00158-021-03025-8}
  {\path{doi:10.1007/s00158-021-03025-8}}.

\bibitem{Dvurechensky2021}
P.~Dvurechensky, E.~Gorbunov, A.~Gasnikov, An accelerated directional
  derivative method for smooth stochastic convex optimization, European Journal
  of Operational Research 290~(2) (2021) 601--621.
\newblock \href {https://doi.org/10.1016/j.ejor.2020.08.027}
  {\path{doi:10.1016/j.ejor.2020.08.027}}.

\bibitem{boole2009}
G.~Boole, A Treatise on the Calculus of Finite Differences, Cambridge
  University Press, Cambridge, 2009.
\newblock \href {https://doi.org/10.1017/CBO9780511693014}
  {\path{doi:10.1017/CBO9780511693014}}.

\bibitem{Martins2003}
J.~R. R.~A. Martins, P.~Sturdza, J.~J. Alonso, The complex-step derivative
  approximation, ACM Transactions on Mathematical Software 29 (2003) 245--262.
\newblock \href {https://doi.org/10.1145/838250.838251}
  {\path{doi:10.1145/838250.838251}}.

\bibitem{KHAN2003}
I.~R. Khan, R.~Ohba, Taylor series based finite difference approximations of
  higher-degree derivatives, Journal of Computational and Applied Mathematics
  154~(1) (2003) 115--124.
\newblock \href {https://doi.org/10.1016/S0377-0427(02)00816-6}
  {\path{doi:10.1016/S0377-0427(02)00816-6}}.

\bibitem{SOD19781}
G.~A. Sod, A survey of several finite difference methods for systems of
  nonlinear hyperbolic conservation laws, Journal of Computational Physics
  27~(1) (1978) 1--31.
\newblock \href {https://doi.org/10.1016/0021-9991(78)90023-2}
  {\path{doi:10.1016/0021-9991(78)90023-2}}.

\bibitem{Martins2013}
J.~R. R.~A. Martins, J.~T. Hwang, Review and unification of methods for
  computing derivatives of multidisciplinary computational models, AIAA Journal
  51~(11) (2013) 2582--2599.
\newblock \href {https://doi.org/doi:10.2514/1.J052184}
  {\path{doi:doi:10.2514/1.J052184}}.

\bibitem{BALAKRISHNA2021}
S.~Balakrishna, W.~W. Schultz, Finite differences for higher order derivatives
  of low resolution data, Mathematics and Computers in Simulation 190 (2021)
  714--722.
\newblock \href {https://doi.org/10.1016/j.matcom.2021.06.011}
  {\path{doi:10.1016/j.matcom.2021.06.011}}.

\bibitem{Griewank1989}
A.~Griewank, Mathematical programming: Recent developments and applications,
  in: M.~Iri, K.~Tanabe (Eds.), On Automatic Differentiation, Kluwer Academic
  Publishers, Dordrecht, 1989, pp. 83--108.

\bibitem{Bischof1992}
C.~Bischof, A.~Carle, G.~Corliss, A.~Griewank, P.~Hovland,
  \href{https://content.iospress.com/articles/scientific-programming/spr1-1-02}{Adifor–generating
  derivative codes from fortran programs}, Scientific Programming 1 (1992)
  11--29.
\newline\urlprefix\url{https://content.iospress.com/articles/scientific-programming/spr1-1-02}

\bibitem{CHINCHALKAR1994197}
S.~Chinchalkar, The application of automatic differentiation to problems in
  engineering analysis, Computer Methods in Applied Mechanics and Engineering
  118~(1) (1994) 197--207.
\newblock \href {https://doi.org/10.1016/0045-7825(94)90113-9}
  {\path{doi:10.1016/0045-7825(94)90113-9}}.

\bibitem{GAWA2003}
A.~Griewank, A.~Walther, Introduction to automatic differentiation, PAMM 2~(1)
  (2003) 45--49.
\newblock \href {https://doi.org/https://doi.org/10.1002/pamm.200310012}
  {\path{doi:https://doi.org/10.1002/pamm.200310012}}.

\bibitem{Neidinger2010}
R.~D. Neidinger, Introduction to automatic differentiation and matlab
  object-oriented programming, SIAM Review 52~(3) (2010) 545--563.
\newblock \href {https://doi.org/10.1137/080743627}
  {\path{doi:10.1137/080743627}}.

\bibitem{Guo2021}
C.~Guo, D.~Poletti, Scheme for automatic differentiation of complex loss
  functions with applications in quantum physics, Phys. Rev. E 103 (2021)
  013309.
\newblock \href {https://doi.org/10.1103/PhysRevE.103.013309}
  {\path{doi:10.1103/PhysRevE.103.013309}}.

\bibitem{OBERBICHLER2021113817}
T.~Oberbichler, R.~Wüchner, K.-U. Bletzinger, Efficient computation of
  nonlinear isogeometric elements using the adjoint method and algorithmic
  differentiation, Computer Methods in Applied Mechanics and Engineering 381
  (2021) 113817.
\newblock \href {https://doi.org/10.1016/j.cma.2021.113817}
  {\path{doi:10.1016/j.cma.2021.113817}}.

\bibitem{CANTUNAVILA2021104086}
K.~Cantún-Avila, D.~González-Sánchez, S.~Díaz-Infante, F.~Peñuñuri,
  Optimizing functionals using differential evolution, Engineering Applications
  of Artificial Intelligence 97 (2021) 104086.
\newblock \href {https://doi.org/10.1016/j.engappai.2020.104086}
  {\path{doi:10.1016/j.engappai.2020.104086}}.

\bibitem{Juedes1991}
D.~W. Juedes, A taxonomy of automatic differentiation tools, in: Griewank,
  G.~F. Corliss (Eds.), Automatic Differentiation of Algorithms: Theory,
  Implementation, and Application, Society for Industrial and Applied
  Mathematics, Philadelphia, 1991, p. 315–329.

\bibitem{Griewank2008}
A.~Griewank, Evaluating Derivatives, Principles and Techniques of Algorithmic
  Differentiation, Vol.~19, Frontiers in Applied Mathematics, SIAM, Portland,
  2008.

\bibitem{Giftthaler2017}
M.~Giftthaler, M.~Neunert, M.~Stäuble, M.~Frigerio, C.~Semini, J.~Buchli,
  Automatic differentiation of rigid body dynamics for optimal control and
  estimation, Advanced Robotics 31~(22) (2017) 1225--1237.
\newblock \href {https://doi.org/10.1080/01691864.2017.1395361}
  {\path{doi:10.1080/01691864.2017.1395361}}.

\bibitem{RAMOS201919}
A.~Ramos, Automatic differentiation for error analysis of monte carlo data,
  Computer Physics Communications 238 (2019) 19--35.
\newblock \href {https://doi.org/https://doi.org/10.1016/j.cpc.2018.12.020}
  {\path{doi:https://doi.org/10.1016/j.cpc.2018.12.020}}.

\bibitem{Peon2020}
R.~{Pe\'on}-Escalante, F.~C. {Jim\'enez}, M.~E. Soberanis, F.~{Pe\~nu\~nuri},
  Path generation with dwells in the optimum dimensional synthesis of
  stephenson {III} six-bar mechanisms, Mechanism and Machine Theory 144 (2020)
  103650.
\newblock \href {https://doi.org/10.1016/j.mechmachtheory.2019.103650}
  {\path{doi:10.1016/j.mechmachtheory.2019.103650}}.

\bibitem{Bolte2020}
J.~Bolte, E.~Pauwels,
  \href{https://proceedings.neurips.cc/paper/2020/file/7a674153c63cff1ad7f0e261c369ab2c-Paper.pdf}{A
  mathematical model for automatic differentiation in machine learning}, in:
  H.~Larochelle, M.~Ranzato, R.~Hadsell, M.~Balcan, H.~Lin (Eds.), Advances in
  Neural Information Processing Systems, Vol.~33, Curran Associates, Inc.,
  2020, pp. 10809--10819.
\newline\urlprefix\url{https://proceedings.neurips.cc/paper/2020/file/7a674153c63cff1ad7f0e261c369ab2c-Paper.pdf}

\bibitem{NIU2021}
Y.~Niu, Y.~He, F.~Xiang, J.~Zhang, Y.~Wu, W.~Tian, G.~Su, S.~Qiu, Automatic
  differentiation approach for solving one-dimensional flow and heat transfer
  problems, Annals of Nuclear Energy 160 (2021) 108361.
\newblock \href {https://doi.org/10.1016/j.anucene.2021.108361}
  {\path{doi:10.1016/j.anucene.2021.108361}}.

\bibitem{ESTEVEZSCHWARZ2021}
D.~Estévez Schwarz, R.~Lamour, Initdae: Computation of consistent values,
  index determination and diagnosis of singularities of daes using automatic
  differentiation in python, Journal of Computational and Applied Mathematics
  387 (2021) 112486.
\newblock \href {https://doi.org/10.1016/j.cam.2019.112486}
  {\path{doi:10.1016/j.cam.2019.112486}}.

\bibitem{Vigliotti2021}
A.~Vigliotti, Automatic differentiation for solid mechanics, Archives of
  Computational Methods in Engineering 28~(3) (2021) 875--895.
\newblock \href {https://doi.org/10.1007/s11831-019-09396-y}
  {\path{doi:10.1007/s11831-019-09396-y}}.

\bibitem{Fakher2022}
B.~Ponsioen, F.~F. Assaad, P.~Corboz, {Automatic differentiation applied to
  excitations with projected entangled pair states}, SciPost Phys. 12 (2022)
  006.
\newblock \href {https://doi.org/10.21468/SciPostPhys.12.1.006}
  {\path{doi:10.21468/SciPostPhys.12.1.006}}.

\bibitem{Kasim2022}
M.~F. Kasim, S.~Lehtola, S.~M. Vinko, Dqc: A python program package for
  differentiable quantum chemistry, The Journal of Chemical Physics 156~(8)
  (2022) 084801.
\newblock \href {https://doi.org/10.1063/5.0076202}
  {\path{doi:10.1063/5.0076202}}.

\bibitem{10.1007/3-540-28438-9_28}
H.~M. B{\"u}cker, G.~F. Corliss, A bibliography of automatic differentiation,
  in: M.~B{\"u}cker, G.~Corliss, U.~Naumann, P.~Hovland, B.~Norris (Eds.),
  Automatic Differentiation: Applications, Theory, and Implementations,
  Springer Berlin Heidelberg, Berlin, Heidelberg, 2006, pp. 321--322.
\newblock \href {https://doi.org/10.1007/3-540-28438-9_28}
  {\path{doi:10.1007/3-540-28438-9_28}}.

\bibitem{Cheng1994}
H.~H. Cheng, Programming with dual numbers and its applications in mechanisms
  design, Engineering with Computers 10~(4) (1994) 212--229.
\newblock \href {https://doi.org/10.1007/BF01202367}
  {\path{doi:10.1007/BF01202367}}.

\bibitem{Jeffrey2012}
J.~A. Fike, J.~J. Alonso, Automatic differentiation through the use of
  hyper-dual numbers for second derivatives, in: S.~Forth, P.~Hovland,
  E.~Phipps, J.~Utke, A.~Walther (Eds.), Recent Advances in Algorithmic
  Differentiation, Springer Berlin Heidelberg, Berlin, Heidelberg, 2012, pp.
  163--173.

\bibitem{Penunuri2013}
F.~Penunuri, R.~Peon-Escalante, C.~Villanueva, O.~Mendoza, C.~A. Cruz-Villar, A
  dual number approach for numerical calculation of velocity and acceleration
  in the spherical 4r mechanism (2013).
\newblock \href {https://doi.org/10.48550/ARXIV.1301.1409}
  {\path{doi:10.48550/ARXIV.1301.1409}}.

\bibitem{Wenbin2013}
W.~Yu, M.~Blair, {DNAD}, a simple tool for automatic differentiation of
  {F}ortran codes using dual numbers, Computer Physics Communications 184
  (2013) 1446--1452.
\newblock \href {https://doi.org/10.1016/j.cpc.2012.12.025}
  {\path{doi:10.1016/j.cpc.2012.12.025}}.

\bibitem{Mendoza2015}
O.~Mendoza-Trejo, C.~A. Cruz-Villar, R.~Peón-Escalante, M.~Zambrano-Arjona,
  F.~Peñuñuri, Synthesis method for the spherical 4r mechanism with minimum
  center of mass acceleration, Mechanism and Machine Theory 93 (2015) 53--64.
\newblock \href {https://doi.org/10.1016/j.mechmachtheory.2015.04.015}
  {\path{doi:10.1016/j.mechmachtheory.2015.04.015}}.

\bibitem{Mex2015}
L.~Mex, C.~A. Cruz-Villar, F.~{Pe\~nu\~nuri}, Closed-form solutions to
  differential equations via differential evolution, Discrete Dynamics in
  Nature and Society 2015 (2015) 910316.
\newblock \href {https://doi.org/10.1155/2015/910316}
  {\path{doi:10.1155/2015/910316}}.

\bibitem{TANAKA2015}
M.~Tanaka, T.~Sasagawa, R.~Omote, M.~Fujikawa, D.~Balzani, J.~Schröder, A
  highly accurate 1st- and 2nd-order differentiation scheme for hyperelastic
  material models based on hyper-dual numbers, Computer Methods in Applied
  Mechanics and Engineering 283 (2015) 22--45.
\newblock \href {https://doi.org/10.1016/j.cma.2014.08.020}
  {\path{doi:10.1016/j.cma.2014.08.020}}.

\bibitem{Moller1993}
M.~F. {M\"oller}, Exact calculation of the product of the hessian matrix of
  feed-forward network error functions and a vector in 0(n) time, DAIMI Report
  Series 22~(432) (1993).
\newblock \href {https://doi.org/10.7146/dpb.v22i432.6748}
  {\path{doi:10.7146/dpb.v22i432.6748}}.

\bibitem{Pearlmutter1994}
B.~A. Pearlmutter, {Fast Exact Multiplication by the Hessian}, Neural
  Computation 6~(1) (1994) 147--160.
\newblock \href {https://doi.org/10.1162/neco.1994.6.1.147}
  {\path{doi:10.1162/neco.1994.6.1.147}}.

\bibitem{Hicken2014}
J.~E. Hicken, {Inexact Hessian-vector products in reduced-space
  differential-equation constrained optimization}, Optimization and Engineering
  15~(3) (2014) 1573--2924.
\newblock \href {https://doi.org/10.1007/s11081-014-9258-6}
  {\path{doi:10.1007/s11081-014-9258-6}}.

\bibitem{Doug2018}
D.~Shi-Dong, S.~Nadarajah, Approximate hessian for accelerated convergence of
  aerodynamic shape optimization problems in an adjoint-based framework,
  Computers \& Fluids 168 (2018) 265--284.
\newblock \href {https://doi.org/10.1016/j.compfluid.2018.04.019}
  {\path{doi:10.1016/j.compfluid.2018.04.019}}.

\bibitem{Song2022}
L.~Song, L.~N. Vicente, {Modeling Hessian-vector products in nonlinear
  optimization: new Hessian-free methods}, IMA Journal of Numerical Analysis
  42~(2) (2021) 1766--1788.
\newblock \href {https://doi.org/10.1093/imanum/drab022}
  {\path{doi:10.1093/imanum/drab022}}.

\bibitem{Schot1978}
S.~H. Schot, Jerk: The time rate of change of acceleration, American Journal of
  Physics 46~(11) (1978) 1090--1094.
\newblock \href {https://doi.org/10.1119/1.11504} {\path{doi:10.1119/1.11504}}.

\bibitem{Eager2016}
D.~Eager, A.-M. Pendrill, N.~Reistad, Beyond velocity and acceleration: jerk,
  snap and higher derivatives, European Journal of Physics 37~(6) (2016)
  065008.
\newblock \href {https://doi.org/10.1088/0143-0807/37/6/065008}
  {\path{doi:10.1088/0143-0807/37/6/065008}}.

\bibitem{Figliolini2019}
G.~Figliolini, C.~Lanni, Jerk and jounce relevance for the kinematic
  performance of long-dwell mechanisms, Advances in Mechanism and Machine
  Science 73 (2016).
\newblock \href {https://doi.org/10.1007/978-3-030-20131-9_22}
  {\path{doi:10.1007/978-3-030-20131-9_22}}.

\bibitem{Fang2020}
Y.~Fang, J.~Qi, J.~Hu, W.~Wang, Y.~Peng, An approach for jerk-continuous
  trajectory generation of robotic manipulators with kinematical constraints,
  Mechanism and Machine Theory 153 (2020) 103957.
\newblock \href {https://doi.org/10.1016/j.mechmachtheory.2020.103957}
  {\path{doi:10.1016/j.mechmachtheory.2020.103957}}.

\bibitem{Hassani1999}
S.~Hassani, Mathematical Physics: A Modern Introduction to Its Foundations,
  Springer-Verlag, New York, 1988.

\bibitem{Carroll2004}
S.~Carroll, Spacetime and Geometry: An Introduction to General Relativity,
  Addison Wesley, New York, 2004.

\bibitem{r_peon2022M}
R.~Pe\'on-Escalante, K.~B. Cant\'un-Avila, O.~Carvente, A.~Espinosa-Romero,
  F.~Pe{\~n}u{\~n}uri, {Dual number implementation to compute higher order
  directional derivatives, {Mendeley} {Data}, v1} (2022).
\newblock \href {https://doi.org/10.17632/kcrm6pmk7d.1}
  {\path{doi:10.17632/kcrm6pmk7d.1}}.

\bibitem{Marban1969}
J.~A. Marban, Directional derivatives in classical optimization, Ph.D. thesis,
  University of Florida (1969).

\bibitem{Hiriart1984}
J.-B. Hiriart-Urruty, Calculus rules on the approximate second-order
  directional derivative of a convex function, SIAM Journal on Control and
  Optimization 22~(3) (1984) 381--404.
\newblock \href {https://doi.org/10.1137/0322025} {\path{doi:10.1137/0322025}}.

\bibitem{Seeger1988}
A.~Seeger, Second order directional derivatives in parametric optimization
  problems, Mathematics of Operations Research 13~(1) (1988) 124--139.
\newblock \href {https://doi.org/10.1287/moor.13.1.124}
  {\path{doi:10.1287/moor.13.1.124}}.

\bibitem{Hassan2012}
H.~Hassan, A.~Mohamad, G.~Atteia, An algorithm for the finite difference
  approximation of derivatives with arbitrary degree and order of accuracy,
  Journal of Computational and Applied Mathematics 236~(10) (2012) 2622--2631.
\newblock \href {https://doi.org/10.1016/j.cam.2011.12.019}
  {\path{doi:10.1016/j.cam.2011.12.019}}.

\bibitem{Clif1873}
W.~Clifford, Preliminary sketch of biquaternions, Proc. London Mathematical
  Society 1~(1-4) (1873) 381--395.
\newblock \href {https://doi.org/10.1112/plms/s1-4.1.381}
  {\path{doi:10.1112/plms/s1-4.1.381}}.

\bibitem{Kalos2021}
L.~Szirmay-Kalos, Higher order automatic differentiation with dual numbers,
  Periodica Polytechnica Electrical Engineering and Computer Science 65~(1)
  (2021) 1–10.
\newblock \href {https://doi.org/10.3311/PPee.16341}
  {\path{doi:10.3311/PPee.16341}}.

\bibitem{Penunuri2016}
F.~Peñuñuri, C.~Cab, O.~Carvente, M.~Zambrano-Arjona, J.~Tapia, A study of
  the classical differential evolution control parameters, Swarm and
  Evolutionary Computation 26 (2016) 86--96.
\newblock \href {https://doi.org/10.1016/j.swevo.2015.08.003}
  {\path{doi:10.1016/j.swevo.2015.08.003}}.

\bibitem{Peon2023}
R.~Pe\'on-Escalante, A.~Espinosa-Romero, F.~{Pe\~nu\~nuri}, Higher order
  kinematic formulas and its numerical computation employing dual numbers,
  Mechanics Based Design of Structures and Machines 0~(0) (2023) 1--16.
\newblock \href {https://doi.org/10.1080/15397734.2023.2203220}
  {\path{doi:10.1080/15397734.2023.2203220}}.

\bibitem{Rico1999}
J.~Rico, J.~Gallardo, J.~Duffy, Screw theory and higher order kinematic
  analysis of open serial and closed chains, Mechanism and Machine Theory
  34~(4) (1999) 559--586.
\newblock \href {https://doi.org/10.1016/S0094-114X(98)00029-9}
  {\path{doi:10.1016/S0094-114X(98)00029-9}}.

\bibitem{Gallardo2014}
J.~Gallardo-Alvarado, Hyper-jerk analysis of robot manipulators, Journal of
  Intelligent \& Robotic Systems 74~(3) (2014) 625--641.
\newblock \href {https://doi.org/10.1007/s10846-013-9849-z}
  {\path{doi:10.1007/s10846-013-9849-z}}.

\bibitem{Muller2014}
A.~{M\"uller}, Higher derivatives of the kinematic mapping and some
  applications, Mechanism and Machine Theory 76 (2014) 70--85.
\newblock \href {https://doi.org/10.1016/j.mechmachtheory.2014.01.007}
  {\path{doi:10.1016/j.mechmachtheory.2014.01.007}}.

\bibitem{Stejskal1996}
V.~Stejskal, M.~Valasek, Kinematics and Dynamics of Machinery (Mechanical
  Engineering), CRC Press, Network, 1996.

\bibitem{Custodio2017}
P.~C. López-Custodio, J.~M. Rico, J.~J. Cervantes-Sánchez, G.~I. Pérez-Soto,
  C.~R. Díez-Martínez, Verification of the higher order kinematic analyses
  equations, European Journal of Mechanics - A/Solids 61~(1) (2017) 198--215.
\newblock \href {https://doi.org/10.1016/j.euromechsol.2016.09.010}
  {\path{doi:10.1016/j.euromechsol.2016.09.010}}.

\end{thebibliography}


\end{document}